\documentclass[11pt, notitlepage]{article}   

\usepackage{amsmath,amsthm,amsfonts}   

\usepackage{amscd}
\usepackage{amsfonts}
\usepackage{amssymb}
\usepackage{color}      
\usepackage{epsfig}
\usepackage{graphicx}           

\usepackage{tikz}
\usepackage{tkz-graph}
\usepackage{tkz-berge}

\theoremstyle{plain}
\newtheorem{theorem}{Theorem}[section]
\newtheorem{lemma}[theorem]{Lemma}

\newtheorem{proposition}[theorem]{Proposition}

\newtheorem{remark}[theorem]{Remark}
\newtheorem{conjecture}[theorem]{Conjecture}

\def\finf{\mathop{{\rm I}\kern -.27 em {\rm F}}\nolimits}

\begin{document}

\date{}

\title{A Comparison between the Metric Dimension and\\ Zero Forcing Number of Trees and Unicyclic Graphs}

\author{{\bf Linda Eroh$^1$, \bf Cong X. Kang$^2$, and \bf Eunjeong Yi$^3$}\\
\small$^1$ University of Wisconsin Oshkosh, Oshkosh, WI 54901, USA\\
\small $^{2,3}$Texas A\&M University at Galveston, Galveston, TX 77553, USA\\
$^1${\small\em eroh@uwosh.edu}; $^2${\small\em kangc@tamug.edu}; $^3${\small\em yie@tamug.edu}}

\maketitle

\begin{abstract}

The \emph{metric dimension} $\dim(G)$ of a graph $G$ is the minimum number of vertices such that every vertex of $G$ is uniquely determined by its vector of distances to the chosen vertices. The \emph{zero forcing number} $Z(G)$ of a graph $G$ is the minimum cardinality of a set $S$ of black vertices (whereas vertices in $V(G)\!\setminus\!S$ are colored white) such that $V(G)$ is turned black after finitely many applications of ``the color-change rule": a white vertex is converted black if it is the only white neighbor of a black vertex. We show that $\dim(T) \leq Z(T)$ for a tree $T$, and that $\dim(G) \le Z(G)+1$ if $G$ is a unicyclic graph; along the way, we characterize trees $T$ attaining $\dim(T)=Z(T)$. For a general graph $G$, we introduce the ``cycle rank conjecture". We conclude with a proof of $\dim(T)-2 \leq \dim(T+e) \le \dim(T)+1$ for $e \in E(\overline{T})$.

\end{abstract}


\section{Introduction}

Let $G = (V(G),E(G))$ be a finite, simple, undirected, connected graph of order $|V(G)|=n \ge 2$ and size $|E(G)|$. The \emph{complement} $\overline{G}$ of a graph $G$ is the graph whose vertex set is $V(G)$ and $uv \in E(\overline{G})$ if and only if $uv \not\in E(G)$ for $u, v \in V(G)$. The \emph{degree} $\deg_G(v)$ of a vertex $v \in V(G)$ is the number of edges incident to the vertex $v$ in $G$; an \emph{end-vertex} is a vertex of degree one. The \emph{distance} between two vertices $u, v \in V(G)$, denoted by $d_G(u, v)$, is the length of a shortest path in $G$ between $u$ and $v$; we omit $G$ when ambiguity is not a concern. 

A vertex $x \in V(G)$ \emph{resolves} a pair of vertices $u,v \in V(G)$ if $d(u,x) \neq d(v,x)$. A set of vertices $W \subseteq V(G)$ \emph{resolves} $G$ if every pair of distinct vertices of $G$ is resolved by some vertex in $W$; then $W$ is called a \emph{resolving set} of $G$. For an ordered set $W=\{w_1, w_2, \ldots, w_k\} \subseteq V(G)$ of distinct vertices, the \emph{metric code} (or \emph{code}, for short) of $v \in V(G)$ with respect to $W$ is the $k$-vector $(d(v, w_1), d(v, w_2), \ldots, d(v, w_k))$; we denote it by $code_W(v)$, and we drop $W$ if the meaning is clear in context. The \emph{metric dimension} of $G$, denoted by $\dim(G)$, is the minimum cardinality over all resolving sets of $G$. Slater \cite{Slater} introduced the concept of a resolving set for a connected graph under the term \emph{locating set}. He referred to a minimum resolving set as a \emph{reference set}, and the cardinality of a minimum resolving set as the \emph{location number} of a graph. Independently, Harary and Melter in~\cite{HM} studied these concepts under the term \emph{metric dimension}. Since metric dimension is suggestive of the dimension of a vector space in linear algebra, sometimes a minimum resolving set of $G$ is called a \emph{basis} of $G$. Metric dimension as a  graph parameter has numerous applications, among them are robot navigation \cite{landmarks}, sonar \cite{Slater}, combinatorial optimization \cite{MathZ}, and pharmaceutical chemistry \cite{CEJO}. In \cite{NPcompleteness}, it is noted that determining the metric dimension of a graph is an NP-hard problem. Metric dimension has been heavily studied; for a survey, see \cite{MDsurvey}. For more on metric dimension in graphs, see \cite{problem, cartesian, joc, Gbar, line, distanceregular, PoZh}.

The notion of a zero forcing set, as well  as the associated zero forcing number, of a simple graph was introduced in~\cite{AIM} to bound the minimum rank for numerous families of graphs. Let each vertex of a graph $G=(V(G), E(G))$ be given one of two colors, dubbed ``black" and ``white" by convention. Let $S$ denote the (initial) set of black vertices of $G$. The \emph{color-change rule} converts the color of a vertex from white to black if the white vertex $u_2$ is the only white neighbor of a black vertex $u_1$; we say that $u_1$ forces $u_2$, which we denote by $u_1 \rightarrow u_2$. And a sequence, $u_1 \rightarrow u_2 \rightarrow \cdots \rightarrow u_{i} \rightarrow u_{i+1} \rightarrow \cdots \rightarrow u_t$, obtained through iterative applications of the color-change rule is called a \emph{forcing chain}. Note that, at each step of the color change, there may be two or more vertices capable of forcing the same vertex. The set $S$ is said to be \emph{a zero forcing set} of $G$ if all vertices of $G$ will be turned black after finitely many applications of the color-change rule. The \emph{zero forcing number} of $G$, denoted by $Z(G)$, is the minimum of $|S|$ over all zero forcing sets $S \subseteq V(G)$. It is known that computing the zero forcing number of a graph is an NP-hard problem (see~\cite{Zhard1, Zhard2}). Zero forcing parameter has been heavily studied; for surveys, see \cite{ZFsurvey, ZFsurvey2}. For more on zero forcing parameter in graphs, see \cite{pathcover, min-degree, iteration, Z+e, ZGbar, proptime}.

More recently, the comparative study of graph parameters is becoming -- it appears -- increasingly fashionable; see \cite{base}, \cite{electronic}, and \cite{mathbohe}, for examples. Our work here is inspired, in part, by these comparative studies. It is also inspired by our observation of the coincidence between zero forcing number and metric dimension for some graphs, as well as by the divergence of these two parameters for other graphs. For graph parameters such as the domination number and the total domination number, which are closely related by their definitions, it is not surprising that there should be inequalities between them. However, metric dimension and zero forcing number arise from rather different contexts and bear no prima facie relation to each other; these facts make any relations discovered between the two parameters all the more interesting and potentially significant. 

The metric dimension and the zero forcing number coincide for paths $P_n$, cycles $C_n$, complete graphs $K_n$, complete bi-partite graphs $K_{s,t}$ ($s+t \ge 3$), for examples; they are $1$, $2$, $n-1$, and $s+t-2$, respectively. For the Cartesian product of two paths and the ``comb" (see Remark~\ref{comb}), the zero forcing number can be seen to be arbitrarily larger than the metric dimension. We will show that $\dim(T) \leq Z(T)$ for a tree $T$, and that $\dim(G) \le Z(G)+1$ if $G$ is a unicyclic graph; both bounds are sharp and, along the way, we characterize trees $T$ attaining $\dim(T)=Z(T)$. On the other hand, the bouquet (or amalgamation) of circles shows that the metric dimension may be arbitrarily larger than the zero forcing number (see \cite{iteration} and \cite{bouquetD}). Nonetheless, we pose the following ``cycle rank conjecture": $\dim(G)\leq Z(G)+r(G)$, where $r(G)$, the cycle rank of $G$, is defined as the minimum number of edges to delete from $G$ so that the resulting graph $G'$ contains no cycle. Towards this conjecture, we show that $\dim(G) \le Z(G)+2r(G)$ if $G$ contains no cycle of even length. We conclude this paper with a proof of $\dim(T)-2 \leq \dim(T+e) \le \dim(T)+1$ for $e \in E(\overline{T})$; see the second paragraph of section 3 for why we include a proof to this known result.


\section{Metric Dimension and Zero Forcing Number of a Tree}

We first recall some results obtained in \cite{CEJO}.

\begin{theorem}[\cite{CEJO}] \label{dimthm}
Let $G$ be a connected graph of order $n \ge 2$. Then
\begin{itemize}
\item[(a)] $\dim(G)=1$ if and only if $G=P_n$,
\item[(b)] $\dim(G)=n-1$ if and only if $G=K_n$,
\item[(c)] for $n \ge 4$, $\dim(G)=n-2$ if and only if $G=K_{s,t}$ ($s,t \ge 1$), $G=K_s + \overline{K}_t$ ($s \ge 1, t \ge 2$), or $G=K_s + (K_1 \cup K_t)$ ($s, t \ge 1$); here, $A+B$ denotes the graph obtained from the disjoint union of graphs $A$ and $B$ by joining every vertex of $A$ with every vertex of $B$.
\end{itemize}
\end{theorem}

\begin{theorem}[\cite{min-degree}]\label{mindegree}
For any connected graph of order $n \ge 2$, $Z(G) \ge \delta(G)$, where $\delta(G)$ is the minimum degree of $G$.
\end{theorem}

\begin{proposition} \label{ObsZ}Let $G$ be a connected graph of order $n \ge 2$. Then
\begin{itemize}
\item[(a)] $Z(G)=1$ if and only if $G=P_n$,
\item[(b)] $Z(G)=n-1$ if and only if $G=K_n$.
\end{itemize}
\end{proposition}

\begin{proof}
(a) Noting that $Z(P_n)=1$ (an end-vertex forms a zero-forcing set of a path), we only need to show that $Z(G)=1$ implies $G=P_n$. Let $\{u_1\}$ be a minimum zero-forcing set of $G$. Then $\deg_{G}(u_1)=1$, otherwise $u_1$ would have more than one white neighbor and not be able to force. Suppose $u_1\rightarrow u_2$ (i.e., $u_1$ forces $u_2$ black); then $u_1$ can no longer force, as it has degree one. Either $n=2$ or $u_2\rightarrow u_3$; the latter implies that $u_2$ must have degree two, as $u_3$ must be the only white neighbor and $u_1$ is the only black neighbor of $u_2$ at this point. Now, we apply this argument inductively until all vertices of $G$ are turned black, and we obtain a forcing chain $u_1\rightarrow u_2 \cdots \rightarrow u_i\rightarrow u_{i+1} \cdots \rightarrow u_{n-1} \rightarrow u_n$. Observe that $u_1$ and $u_n$ each has degree one, whereas each $u_i$ for $1<i<n$ has degree two: $u_i$ forcing $u_{i+1}$ means that $u_{i+1}$ is the only white neighbor of $u_i$, which, while forcing, has only one black neighbor $u_{i-1}$. This means that $G$ is $P_n$.

(b) Note that $Z(K_n)=n-1$ (all but one vertex of $G$ forms a zero-forcing set of a complete graph). On the other hand, $Z(G)=n-1$ implies $G=K_n$. For $V(G)=\{u_1, u_2, \ldots, u_{n}\}$, suppose $G\not\cong K_n$ (thus $n>2$) and let $e=u_{1}u_{n}\not\in E(G)$. By relabeling if necessary, we may assume that $u_1u_{n-1} \in E(G)$ and $u_tu_n \in E(G)$, where $2 \le t \le n-1$. Then the set $\{u_1, u_2, \ldots, u_{n-2}\}$ is a zero forcing set of cardinality $n-2$ for $G$: $u_1 \rightarrow u_{n-1}$; subsequent to (if not simultaneous with) $u_{n-1}$ turning black, $u_t \rightarrow u_n$.~\hfill
\end{proof}

\begin{proposition}
Let $S_0$ be a zero-forcing set of a connected graph $G$. If the entire vertex set of $G$ turns black after one global application of the color-change rule, then $S_0$ is a resolving set for $G$.
\end{proposition}

\begin{proof}
Let $S_1 = V(G) \setminus S_0$. Let $S_0 = \{v_1, v_2, \ldots, v_t\}$; notice that each $v_i \in S_0$ ($1 \le i \le t$) has a unique code, since it is the only code with $0$ in the $i$-th coordinate. If $x \in S_1$, then there exists $v_j \in S_0$ such that $x$ is the only neighbor of $v_j$ that is not in $S_0$; thus, $x$ is the only vertex with $1$ in the $j$-th coordinate of its code and no zero in its code, so it has a unique code.~\hfill
\end{proof}

The following definitions are introduced in \cite{CEJO}. Fix a graph $G$. A vertex of degree at least three is called a \emph{major vertex}. An end-vertex $u$ is called \emph{a terminal vertex of a major vertex} $v$ if $d(u, v)<d(u, w)$ for every other major vertex $w$. The \emph{terminal degree} of a major vertex $v$ in $G$, denoted by $ter_G(v)$, is the number of terminal vertices of $v$. A major vertex $v$ is an \emph{exterior major vertex} (emv) if it has positive terminal degree. Let $\sigma(G)$ denote the sum of terminal degrees of all major vertices of $G$, and let $ex(G)$ denote the number of emvs of $G$. We further define an \emph{exterior degree two vertex} to be a vertex of degree 2 that lies on a path from a terminal vertex to its major vertex, and an \emph{interior degree two vertex} to be a vertex of degree 2 such that the shortest path to any terminal vertex includes a major vertex.  We refer to the components of $G-v$ as the \emph{branches} of $v$.  If $v$ is an emv, then a branch which contains a terminal vertex of $v$ will be called an \emph{exterior branch} of $v$. Two vertices $u, v \in V(G)$ are called \emph{twins} if $N(u) \setminus \{v\}=N(v) \setminus \{u\}$, where $N(u)$ is the set of all vertices adjacent to $u$ in $G$; notice that for any set $S$ with $S \cap \{u,v\} = \emptyset$, $code_S(u)=code_S(v)$.

\begin{lemma}[\cite{CEJO}] \label{dim geq sigma-ex}
If $G$ is any graph, then $\dim(G)\geq\sigma(G)-ex(G)$.
\end{lemma}

\begin{theorem}[\cite{CEJO, landmarks, PoZh}] \label{tree}
If $T$ is a tree that is not a path, then $\dim(T)=\sigma(T)-ex(T)$.
\end{theorem}

\begin{theorem} \label{DZtree}
For any tree $T$, $\dim(T) \leq Z(T)$.
\end{theorem}

\begin{proof}
If $T$ is a path, $\dim(T)=Z(T)=1$ by (a) of Theorem~\ref{dimthm} and (a) of Proposition~\ref{ObsZ}. So, we only need to consider trees that are not paths. Take any set $S \subseteq V(T)$ with $|S| < \sigma(T)-ex(T)$. There then must be an emv $u$ and a pair of terminal vertices $x$ and $y$ of $u$ such that no vertex on the path from $u$ to $x$ (except possibly $u$) is in $S$, and no vertex on the path from $u$ to $y$ (except possibly $u$) is in $S$. Let $x'$ be the vertex adjacent to $u$ on the $u-x$ path and $y'$ be the vertex adjacent to $u$ on the $u-y$ path. Consider iterative applications of the color-change rule to the initial black set $S$. Notice that even if $u$ is turned black, both $x'$ and $y'$ will remain white; so $S$ cannot be a zero-forcing set. Thus, any zero-forcing set of $T$ must have cardinality at least $\sigma(T)-ex(T)=\dim(T)$.~\hfill
\end{proof}

The \emph{path cover number} $P(G)$ of $G$ is the minimum number of vertex disjoint paths, occurring as induced subgraphs of $G$, that cover all the vertices of $G$.

\begin{theorem}[\cite{AIM, pathcover}]\label{pathcover}
\begin{itemize}
\item[(a)] \cite{pathcover} For any graph $G$, $P(G) \le Z(G)$.
\item[(b)] \cite{AIM} For any tree $T$, $P(T) = Z(T)$.
\end{itemize}
\end{theorem}

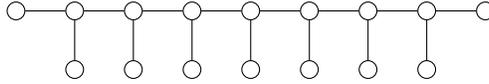
\begin{figure}[ht]
\centering
\begin{tikzpicture}[scale=.65, transform shape]
\node [draw, shape=circle] (1) at  (-1,1) {};
\node [draw, shape=circle] (2) at  (0.2,1) {};
\node [draw, shape=circle] (3) at  (1.4,1) {};
\node [draw, shape=circle] (4) at  (2.6,1) {};
\node [draw, shape=circle] (5) at  (3.8,1) {};
\node [draw, shape=circle] (6) at  (5,1) {};
\node [draw, shape=circle] (7) at  (6.2,1) {};
\node [draw, shape=circle] (8) at  (7.4,1) {};
\node [draw, shape=circle] (9) at  (8.6,1) {};

\node [draw, shape=circle] (22) at  (0.2,-0.2) {};
\node [draw, shape=circle] (33) at  (1.4,-0.2) {};
\node [draw, shape=circle] (44) at  (2.6,-0.2) {};
\node [draw, shape=circle] (55) at  (3.8,-0.2) {};
\node [draw, shape=circle] (66) at  (5,-0.2) {};
\node [draw, shape=circle] (77) at  (6.2,-0.2) {};
\node [draw, shape=circle] (88) at  (7.4,-0.2) {};

\draw(1)--(2)--(3)--(4)--(5)--(6)--(7)--(8)--(9);
\draw(2)--(22);
\draw(3)--(33);
\draw(4)--(44);
\draw(5)--(55);
\draw(6)--(66);
\draw(7)--(77);
\draw(8)--(88);

\end{tikzpicture}
\caption{A tree $T$ showing that $Z(T)$ can be arbitrarily bigger than $\dim(T)$}\label{fig1}
\end{figure}

\begin{remark}\label{comb}
We note that $Z(G)$ can be arbitrarily larger than $\dim(G)$ for a graph $G$. It's shown in~\cite{cartesian} that the metric dimension of the Cartesian product of paths $P_m \Box P_n$ is two, whereas $Z(P_m \Box P_n)=\min\{m,n\}$, as shown in~\cite{AIM}. The tree $T$ in Figure~\ref{fig1} is another example. By Theorem~\ref{tree}, $\dim(T)=2$. On the other hand, $Z(T)=P(T)=5$, since each path contains at most two end-vertices. Clearly, by adding more $P_2$'s to the horizontal path, we can arbitrarily increase the zero forcing number while holding the metric dimension at two.
\end{remark}

Next, we characterize trees $T$ satisfying $\dim(T)=Z(T)$. 

\begin{theorem}\label{dim=Z,T}
For any tree $T$, we have $\dim(T)=Z(T)$ if and only if $T$ has no interior degree two vertices and each major vertex $v$ of $T$ satisfies $ter_T(v) \ge 2$.
\end{theorem}

\begin{proof}
($\Longleftarrow$) If $ex(T)=0$, then $T$ is a path, and we have $\dim(T)=Z(T)=1$. So, we consider $ex(T) \ge 1$. Let $v_1, v_2, \ldots, v_m$ be all the emvs of $T$ with $ter_T(v_i) \ge 2$ for each $i$, $1 \le i \le m$; further, suppose $T$ has no interior degree 2 vertices. Denote by $\ell_{i,1}, \ell_{i,2}, \ldots, \ell_{i, k_i}$ the terminal vertices of $v_i$. There exists a path $\ell_{i,1}, \ldots, v_i, \ldots, \ell_{i,2}$ for each $v_i$, since $k_i \ge 2$ for each $i$; there are also $(k_i-2)$ additional paths associated to each $v_i$ (noting that $v_i$ may belong to only one path). So, $P(T) \le \sum_{i=1}^{m}(k_i-1)=(\sum_{i=1}^{m}k_i)-m=\sigma(T)-ex(T)=\dim(T)$. On the other hand, $\dim(T) \le Z(T)=P(T)$ by Theorem \ref{DZtree} and (b) of Theorem \ref{pathcover}. Thus, $\dim(T)=Z(T)$.

($\Longrightarrow$) Let $\dim(T)=Z(T)$. A path trivially satisfies the conditions on $T$; so, let $T$ be a tree which is not a path. Suppose $T$ has an interior degree two vertex $u_0$ or a major vertex $v_0$ with $ter_T(v_0)<2$. Since $Z(T)=P(T)$, the desired contradiction is reached if we exhibit a subcover $B$ of a minimum path cover of $T$ such that $B\cap\{u_0,v_0\}=\emptyset$ and $|B|=\dim(T)=\sigma(T)-ex(T)$. To this end, let $v_1, \ldots, v_k$ enumerate all major vertices with terminal degree at least two; notice $k\geq 1$ as $T$ is not a path. Put $ter_T(v_i)=m_i$, and let $\ell_{i_j},\ldots, w_{i_j}, v_i$ (where $1\leq j\leq m_i$) denote the path between $v_i$ and its $j$-th terminal vertex (a leaf) $\ell_{i_j}$. Let $B_i=\{\{\ell_{i_1},\ldots,w_{i_1},v_i,w_{i_2},\ldots,\ell_{i_2}\},\{\ell_{i_3},\ldots,w_{i_3}\},\cdots,\{\ell_{i_{m_i}},\ldots,w_{i_{m_i}}\}\}$. $B_i$ belongs to a minimum path cover, since $|B_i|=m_i-1$ and the use of a path containing $v_i$ but fewer than two terminal vertices of $v_i$ necessitates the use of at least $m_i-1$ paths to cover the remaining vertices covered by $B_i$. It is then clear that $B=\bigcup_{1\leq i\leq k}B_i$ is the requisite subcover.~\hfill
\end{proof}


\section{Metric Dimension and Zero Forcing Number of a Unicyclic Graph}

The \emph{cycle rank} of a graph $G$, denoted by $r(G)$, is defined as $|E(G)|-|V(G)|+1$. For a tree $T$, $r(T)=0$. If a graph $G$ has $r(G)=1$, we call it a \emph{unicyclic} graph. By $T+e$, we shall mean a unicyclic graph obtained from a tree $T$ by attaching a new edge $e \in E(\overline{T})$. In~\cite{PoZh}, the notion of a resolving set $W$ with the property $code_W(u)-code_W(v) \neq (a, \dots, a)$ for any $a\in \mathbb{Z}$ was identified and shown to be very useful. We will say that ``$G$ is \emph{strongly resolved} by $W$" if  $code_W(u)-code_W(v) \neq (a, \dots, a)$ for any $a\in \mathbb{Z}$ and any $u, v\in V(G)$. Still following~\cite{PoZh}, observe that $u\!\sim_{W}\!v$ if and only if $code_W(u)-code_W(v)=(a, \dots, a)$ for some $a\in \mathbb{Z}$ defines an equivalence relation $\sim_W$ on $V(G)$; let $[u]_W$ denote the equivalence class of $u$ under this relation.

The upper bound in the following theorem (Theorem~\ref{unicyclic}) is fundamental to Theorem~\ref{bigtheorem}. It is stated in~\cite{CEJO}, along with an outline of proof; unfortunately, the outline is logically flawed (see Remarks~\ref{ce1CEJO} and~\ref{ce2CEJO}). The theorem is also attributed to~\cite{PoZh} by some authors, but we do not see it as an immediate (and unstated) corollary of~\cite{PoZh}. In consideration of these facts, we will include a proof to Theorem~\ref{unicyclic} in the final section of this paper for self-containedness. Our proof to the lower bound in Theorem~\ref{unicyclic} is a modification of that given in~\cite{CEJO}, and our proof to the upper bound in Theorem~\ref{unicyclic} is based on some of the ideas contained in~\cite{PoZh}.

\begin{theorem}\label{unicyclic}
If $T$ is a tree of order at least three and $e$ is an edge of $\overline{T}$, then $$\dim(T)-2 \le \dim(T+e) \le \dim(T)+1.$$
\end{theorem}

\begin{theorem}[\cite{Z+e}]\label{Z+}
Let G be any graph. 
\begin{itemize}
\item[(a)] For $v \in V(G)$, $Z(G)-1 \le Z(G-\{v\}) \le Z(G)+1$.
\item[(b)] For $e \in E(G)$, $Z(G)-1 \le Z(G-e) \le Z(G)+1$.
\end{itemize}
\end{theorem}

\begin{remark}
As an immediate consequence of Theorem~\ref{DZtree}, Theorem~\ref{unicyclic}, and (b) of Theorem~\ref{Z+}, we have $\dim(T+e) \le Z(T+e)+2$, where $e \in E(\overline{T})$. In order for $\dim(T+e)=Z(T+e)+2$, $T$ must satisfy $\dim(T)=Z(T)$; further, $T+e$ must satisfy both $\dim(T+e)=\dim(T)+1$ and $Z(T+e)=Z(T)-1$: we will show that this can not happen.
\end{remark}

\begin{remark}\label{ex00}
A tree $T$ with $ex(T)=0$ is a path $P_n$. One easily sees that $Z(P_n+e)=2$, and it follows from Theorem 4.2 of~\cite{PoZh} that $\dim(P_n+e)\leq 3$. In fact, $\dim(P_n+e)=2$, since the two end-vertices of $P_n$ always form a resolving set for $P_n+e$. We will prove Theorem~\ref{bigtheorem} by inducting on $ex(T)$ of a tree $T$ satisfying $\dim(T)=Z(T)$. To facilitate the induction process, we will first establish the result when $ex(T)=1$.
\end{remark}

\begin{proposition}\label{ex1}
Let $T$ be a tree with $ex(T)=1$. Then $\dim(T+e) \le Z(T+e)+1$, where $e \in E(\overline{T})$.
\end{proposition}

\begin{proof}
If $\dim(T)<Z(T)$, then $\dim(T+e) \le Z(T+e)+1$ by Theorem~\ref{unicyclic} and by (b) of Theorem~\ref{Z+}. So, by Theorem~\ref{DZtree}, we only need to consider $T$ satisfying $\dim(T)=Z(T)$. By Theorem~\ref{dim=Z,T} and by the condition $ex(T)=1$, if $v$ is the unique (exterior) major vertex of $T$, then the terminal degree of $v$ is at least three. Let $\ell_1, \ell_2, \ldots, \ell_j$ be the terminal vertices of $v$ in $T$, where $j \ge 3$. Let $B=\{x_2,\ldots, x_j\}$ where, for $2\leq i\leq j$, each $x_i$ is a vertex lying on the $\ell_i-v$ path and not equal to $v$. ($B$ is, if you will, a set of functions which specializes to the prescribed set of vertices when (partially) included in a zero forcing set.) Notice that $\dim(T)=Z(T)=j-1$. Let $C$ be the unique cycle in $T+e$. Let $S$ be a zero-forcing set of $T+e$, and we consider two cases.

\emph{Case 1. $C$ does not contain $v$:}  Let $s$ and $s'$ be degree 2 vertices that lie on the path from $\ell_1$ to $v$. Without loss of generality (WLOG), we need to consider two cases: (A) $e=ss'$ (see (a) of Figure~\ref{fig2}) and (B) $e=s\ell_1$ (see (b) of Figure~\ref{fig2}). In each case, $S$ must contain all but one element of the set $B$; WLOG, let $S\supseteq S_0=\{\ell_2, \ell_3, \ldots, \ell_{j-1}\}$. If $S_0$ is the initial black set of $T+e$, once the vertex $v$ is turned black, $v$ has two white neighbors and can not force. Thus, at least a vertex in $V(T) \setminus S_0$ must belong to $S$, and hence $Z(T+e) \ge j-1$. Since $Z(T+e) \ge Z(T)$, we have $\dim(T+e) \le Z(T+e)+1$ by Theorem~\ref{unicyclic}.

\begin{figure}[ht]
\centering
\begin{tikzpicture}[scale=.65, transform shape]
\node [draw, shape=circle] (a0) at  (0,2) {};
\node [draw, shape=circle] (a1) at  (-2,1) {};
\node [draw, shape=circle] (a2) at  (-2,0) {};
\node [draw, shape=circle] (a3) at  (-2,-1) {};
\node [draw, shape=circle] (a4) at  (-2,-2.3) {};
\node [draw, shape=circle] (a11) at  (-0.5,1) {};
\node [draw, shape=circle] (a22) at  (-0.5,0) {};
\node [draw, shape=circle] (a111) at  (0.5,1) {};
\node [draw, shape=circle] (a222) at  (0.5,0) {};
\node [draw, shape=circle] (a1111) at  (2,1) {};

\draw(a1)--(a2)--(a3);
\draw(a0)--(a11)--(a22);
\draw(a222)--(a111)--(a0)--(a1111);
\draw(a1)--(-1.55, 1.24);
\draw(-0.5,1.74)--(a0);
\draw[thick, style=dotted](-1.55,1.24)--(-0.5, 1.74);
\draw(a3)--(-2, -1.4);
\draw(-2,-1.9)--(a4);
\draw[thick, style=dotted] (-2,-1.4)--(-2,-1.9); 
\draw[thick, style=dotted] (0.9,1)--(1.6,1);
\draw[very thick](a1) .. controls (-1.3,0.3) and (-1.3, -0.3) .. (a3);

\node [scale=1.4] at (0,2.4) {$v$};
\node [scale=1.4] at (-2.5,1) {$s$};
\node [scale=1.4] at (-2.5,-0.9) {$s'$};
\node [scale=1.4] at (-2.5,-2.3) {$\ell_1$};
\node [scale=1.4] at (-0.4,-0.5) {$\ell_2$};
\node [scale=1.4] at (0.6,-0.5) {$\ell_3$};
\node [scale=1.4] at (2.1,0.5) {$\ell_j$};
\node [scale=1.4] at (0,-3) {\large (a) $T+e$};

\node [draw, shape=circle] (b0) at  (8,2) {};
\node [draw, shape=circle] (b1) at  (6,1) {};
\node [draw, shape=circle] (b2) at  (6,0) {};
\node [draw, shape=circle] (b3) at  (6,-1) {};
\node [draw, shape=circle] (b4) at  (6,-2.3) {};
\node [draw, shape=circle] (b11) at  (7.5,1) {};
\node [draw, shape=circle] (b22) at  (7.5,0) {};
\node [draw, shape=circle] (b111) at  (8.5,1) {};
\node [draw, shape=circle] (b222) at  (8.5,0) {};
\node [draw, shape=circle] (b1111) at  (10,1) {};

\draw(b1)--(b2)--(b3);
\draw(b0)--(b11)--(b22);
\draw(b222)--(b111)--(b0)--(b1111);
\draw(b1)--(6.45, 1.24);
\draw(7.5,1.74)--(b0);
\draw[thick, style=dotted](6.45,1.24)--(7.5, 1.74);
\draw(b3)--(6, -1.4);
\draw(6,-1.9)--(b4);
\draw[thick, style=dotted] (6,-1.4)--(6,-1.9); 
\draw[thick, style=dotted] (8.9,1)--(9.6,1);
\draw[very thick](b1) .. controls (6.7,0) and (6.7, -1.3) .. (b4);

\node [scale=1.4] at (8,2.4) {$v$};
\node [scale=1.4] at (5.5,1) {$s$};
\node [scale=1.4] at (5.5,-2.3) {$\ell_1$};
\node [scale=1.4] at (7.6,-0.5) {$\ell_2$};
\node [scale=1.4] at (8.6,-0.5) {$\ell_3$};
\node [scale=1.4] at (10.1,0.5) {$\ell_j$};
\node [scale=1.4] at (8,-3) {\large (b) $T+e$};
\end{tikzpicture}
\caption{Unicyclic graphs $T+e$ with $ex(T)=1$ and $v \not\in V(C)$}\label{fig2}
\end{figure}

\emph{Case 2. $C$ contains $v$:} Let $s_1$ ($s_2$, respectively) be an exterior degree two vertex that lies on the path from $\ell_1$ to $v$ ($\ell_2$ to $v$, respectively). WLOG, we need to consider five cases: (A) $e=\ell_1v$ (see (a) of Figure~\ref{fig3}); (B) $e=s_1v$ (see (b) of Figure~\ref{fig3}); (C) $e=\ell_1\ell_2$ (see (c) of Figure~\ref{fig3}); (D) $e=s_1\ell_2$ (see (d) of Figure~\ref{fig3}); (E) $e=s_1s_2$ (see (e) of Figure~\ref{fig3}).

First, we consider (A) and (B). Note that $S$ must contain all but one element of the set $B$; WLOG, let $S\supseteq S_0=\{\ell_2, \ell_3, \ldots, \ell_{j-1}\}$. If $S_0$ is the initial black set for $T+e$, once $v$ is turned black, $v$ has three white neighbors and can not force. So, at least a vertex in $V(T) \setminus S_0$ must belong to $S$, and hence $Z(T+e) \ge Z(T)=j-1$. Thus, $\dim(T+e) \le Z(T+e)+1$ by Theorem~\ref{unicyclic}.

Second, we consider (C), (D), and (E). If $j=3$, then it's obvious that $Z(T+e)\geq 2=Z(T)$, and we have $\dim(T+e) \le Z(T+e)+1$. So, we consider $j \ge 4$. Let $W=\{u_1, u_2, \ell_3, \ldots, \ell_{j-1}\}$, where $u_1$ and $u_2$ are vertices lying on the unique cycle $C$ such that $d_{T+e}(u_1,u_2)$ is the diameter of $C$ and $v\not\in \{u_1,u_2\}$. We contend that $W$ is a resolving set for $T+e$, and thus $\dim(T+e)\leq \dim(T)$. Since $j\geq 4$, $|W|\geq 3$. Note, as in~\cite{PoZh}, that vertices on $C$ are ``strongly resolved" by $\{u_1,u_2,v\}$ (no vertex on $C$ is simultaneously closer to all three vertices, \emph{as chosen}, than another vertex on $C$). Hence, for $\ell \in \{\ell_3, \ldots, \ell_{j-1}\}$, the vertices on $C$ are also strongly resolved by $\{u_1,u_2,\ell \}$ such that $d_{T+e}(u,\ell)=d_{T+e}(u,v)+d_{T+e}(v,\ell)$, whenever $u$ is a vertex lying on $C$. Viewing $T+e$ as a collection of trees $T_i$ rooted at vertices of $C$ (together with $C$), the strong resolution of $C$ by $\{u_1, u_2, v\}$ ensures that no vertex on $T_i$ will have the same code as a vertex on $T_j$ if $i\neq j$. Further, vertices on the $j-2$ paths rooted at $v$ are clearly resolved among themselves. Thus, $\dim(T+e) \le \dim(T)$, and hence $\dim(T+e) \le Z(T+e)+1$ by (b) of Theorem~\ref{Z+}.~\hfill
\end{proof}

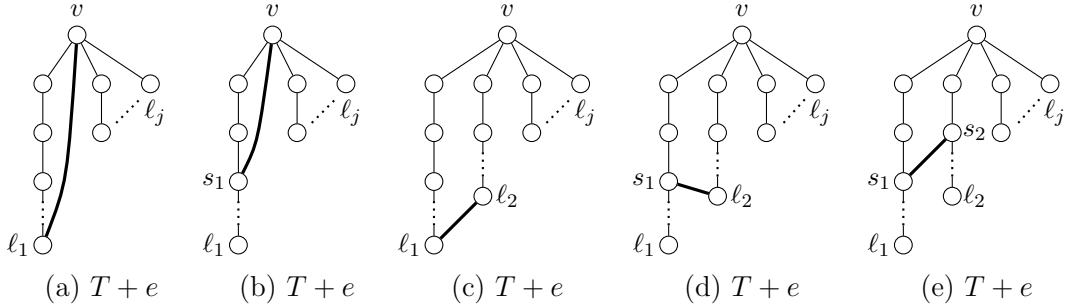
\begin{figure}[ht]
\centering
\begin{tikzpicture}[scale=.65, transform shape]
\node [draw, shape=circle] (a0) at  (0,3) {};
\node [draw, shape=circle] (a1) at  (-0.7,2) {};
\node [draw, shape=circle] (a2) at  (-0.7,1) {};
\node [draw, shape=circle] (a3) at  (-0.7,0) {};
\node [draw, shape=circle] (a4) at  (-0.7,-1.3) {};
\node [draw, shape=circle] (a11) at  (0.5,2) {};
\node [draw, shape=circle] (a22) at  (0.5,1) {};
\node [draw, shape=circle] (a111) at  (1.5,2) {};

\draw(a0)--(a1)--(a2)--(a3);
\draw(a22)--(a11)--(a0)--(a111);
\draw(a3)--(-0.7, -0.4);
\draw(-0.7,-0.9)--(a4);
\draw[thick, style=dotted](-0.7,-0.4)--(-0.7,-0.9);
\draw[thick, style=dotted] (0.8,1.2)--(1.3,1.7);
\draw[very thick](a4) .. controls (-0.2,0) .. (a0);

\node [scale=1.4] at (0,3.5) {$v$};
\node [scale=1.4] at (-1.2,-1.3) {$\ell_1$};
\node [scale=1.4] at (1.6,1.4) {$\ell_j$};
\node [scale=1.4] at (0.5,-2.2) {\large (a) $T+e$};

\node [draw, shape=circle] (b0) at  (4,3) {};
\node [draw, shape=circle] (b1) at  (3.3,2) {};
\node [draw, shape=circle] (b2) at  (3.3,1) {};
\node [draw, shape=circle] (b3) at  (3.3,0) {};
\node [draw, shape=circle] (b4) at  (3.3,-1.3) {};
\node [draw, shape=circle] (b11) at  (4.5,2) {};
\node [draw, shape=circle] (b22) at  (4.5,1) {};
\node [draw, shape=circle] (b111) at  (5.5,2) {};

\draw(b0)--(b1)--(b2)--(b3);
\draw(b22)--(b11)--(b0)--(b111);
\draw(b3)--(3.3, -0.4);
\draw(3.3,-0.9)--(b4);
\draw[thick, style=dotted](3.3,-0.4)--(3.3,-0.9);
\draw[thick, style=dotted] (4.8,1.2)--(5.3,1.7);
\draw[very thick](b3) .. controls (3.7,0.8) .. (b0);

\node [scale=1.4] at (4,3.5) {$v$};
\node [scale=1.4] at (2.8,0) {$s_1$};
\node [scale=1.4] at (2.8,-1.3) {$\ell_1$};
\node [scale=1.4] at (5.6,1.4) {$\ell_j$};
\node [scale=1.4] at (4.5,-2.2) {\large (b) $T+e$};

\node [draw, shape=circle] (c0) at  (8.8,3) {};
\node [draw, shape=circle] (c1) at  (7.3,2) {};
\node [draw, shape=circle] (c2) at  (7.3,1) {};
\node [draw, shape=circle] (c3) at  (7.3,0) {};
\node [draw, shape=circle] (c4) at  (7.3,-1.3) {};
\node [draw, shape=circle] (c11) at  (8.3,2) {};
\node [draw, shape=circle] (c22) at  (8.3,1) {};
\node [draw, shape=circle] (c33) at  (8.3,-0.3) {};
\node [draw, shape=circle] (c111) at  (9.3,2) {};
\node [draw, shape=circle] (c222) at  (9.3,1) {};
\node [draw, shape=circle] (c1111) at  (10.3,2) {};

\draw(c3)--(c2)--(c1)--(c0)--(c11)--(c22);
\draw(c222)--(c111)--(c0)--(c1111);
\draw(c3)--(7.3, -0.4);
\draw(7.3,-0.9)--(c4);
\draw[thick, style=dotted](7.3,-0.4)--(7.3,-0.9);
\draw(c22)--(8.3, 0.6);
\draw(8.3,0.1)--(c33);
\draw[thick, style=dotted](8.3,0.6)--(8.3,0.1);
\draw[thick, style=dotted] (9.6,1.2)--(10.1,1.7);
\draw[very thick](c4)--(c33);

\node [scale=1.4] at (8.8,3.5) {$v$};
\node [scale=1.4] at (6.8,-1.3) {$\ell_1$};
\node [scale=1.4] at (8.8,-0.3) {$\ell_2$};
\node [scale=1.4] at (10.4,1.4) {$\ell_j$};
\node [scale=1.4] at (8.8,-2.2) {\large (c) $T+e$};

\node [draw, shape=circle] (d0) at  (13.6,3) {};
\node [draw, shape=circle] (d1) at  (12.1,2) {};
\node [draw, shape=circle] (d2) at  (12.1,1) {};
\node [draw, shape=circle] (d3) at  (12.1,0) {};
\node [draw, shape=circle] (d4) at  (12.1,-1.3) {};
\node [draw, shape=circle] (d11) at  (13.1,2) {};
\node [draw, shape=circle] (d22) at  (13.1,1) {};
\node [draw, shape=circle] (d33) at  (13.1,-0.3) {};
\node [draw, shape=circle] (d111) at  (14.1,2) {};
\node [draw, shape=circle] (d222) at  (14.1,1) {};
\node [draw, shape=circle] (d1111) at  (15.1,2) {};

\draw(d3)--(d2)--(d1)--(d0)--(d11)--(d22);
\draw(d222)--(d111)--(d0)--(d1111);
\draw(d3)--(12.1, -0.4);
\draw(12.1,-0.9)--(d4);
\draw[thick, style=dotted](12.1,-0.4)--(12.1,-0.9);
\draw(d22)--(13.1, 0.6);
\draw(13.1,0.1)--(d33);
\draw[thick, style=dotted](13.1,0.6)--(13.1,0.1);
\draw[thick, style=dotted] (14.45,1.2)--(14.9,1.7);
\draw[very thick](d3)--(d33);

\node [scale=1.4] at (13.6,3.5) {$v$};
\node [scale=1.4] at (11.6,-1.3) {$\ell_1$};
\node [scale=1.4] at (11.6,0) {$s_1$};
\node [scale=1.4] at (13.6,-0.3) {$\ell_2$};
\node [scale=1.4] at (15.2,1.4) {$\ell_j$};
\node [scale=1.4] at (13.6,-2.2) {\large (d) $T+e$};

\node [draw, shape=circle] (e0) at  (18.4,3) {};
\node [draw, shape=circle] (e1) at  (16.9,2) {};
\node [draw, shape=circle] (e2) at  (16.9,1) {};
\node [draw, shape=circle] (e3) at  (16.9,0) {};
\node [draw, shape=circle] (e4) at  (16.9,-1.3) {};
\node [draw, shape=circle] (e11) at  (17.9,2) {};
\node [draw, shape=circle] (e22) at  (17.9,1) {};
\node [draw, shape=circle] (e33) at  (17.9,-0.3) {};
\node [draw, shape=circle] (e111) at  (18.9,2) {};
\node [draw, shape=circle] (e222) at  (18.9,1) {};
\node [draw, shape=circle] (e1111) at  (19.9,2) {};

\draw(e3)--(e2)--(e1)--(e0)--(e11)--(e22);
\draw(e222)--(e111)--(e0)--(e1111);
\draw(e3)--(16.9, -0.4);
\draw(16.9,-0.9)--(e4);
\draw[thick, style=dotted](16.9,-0.4)--(16.9,-0.9);
\draw(e22)--(17.9, 0.6);
\draw(17.9,0.1)--(e33);
\draw[thick, style=dotted](17.9,0.6)--(17.9,0.1);
\draw[thick, style=dotted] (19.25,1.2)--(19.7,1.7);
\draw[very thick](e3)--(e22);

\node [scale=1.4] at (18.4,3.5) {$v$};
\node [scale=1.4] at (16.4,-1.3) {$\ell_1$};
\node [scale=1.4] at (16.4,0) {$s_1$};
\node [scale=1.4] at (18.36,-0.3) {$\ell_2$};
\node [scale=1.4] at (18.36,1) {$s_2$};
\node [scale=1.4] at (20,1.4) {$\ell_j$};
\node [scale=1.4] at (18.4,-2.2) {\large (e) $T+e$};

\end{tikzpicture}
\caption{Unicyclic graphs $T+e$ with $ex(T)=1$ and $v \in V(C)$}\label{fig3}
\end{figure}

\begin{theorem}\label{bigtheorem}
If $T$ is a tree and $e \in E(\overline{T})$, then
\begin{equation}\label{D,Z+1}
\dim(T+e) \le Z(T+e)+1.
\end{equation}
\end{theorem}

\begin{proof}
If $\dim(T)<Z(T)$, then $\dim(T+e) \le Z(T+e)+1$ by Theorem~\ref{unicyclic} and by (b) of Theorem~\ref{Z+}. So, by Theorem~\ref{DZtree}, we only need to consider trees $T$ satisfying $\dim(T)=Z(T)$; we will induct on $ex(T)$ for such $T$'s. We have already established inequality~(\ref{D,Z+1}) when $ex(T)\leq 1$ with Remark~\ref{ex00} and Proposition~\ref{ex1}. So, assume that  (\ref{D,Z+1}) holds for any tree $T$ with $ex(T)=k \ge 1$, and consider a tree $T$ with $ex(T)=k+1 \ge 2$. By Theorem~\ref{dim=Z,T}, we need only to consider trees $T$ that has no interior degree 2 vertex and no major vertex with terminal degree less than 2. The key idea in this proof is rewriting $T+e$ as $T'+e'$, where $e'$ is some edge on the unique cycle of $T+e$ such that, typically, either $T'$ satisfies $ex(T')<ex(T)$ (then induction hypothesis applies) or $T'$ contains a structural element such as a vertex of interior degree 2 or a major vertex of terminal degree less than 2 (thus $\dim(T')<Z(T')$, and then $\dim(T+e)\leq Z(T+e)+1$, since $T'+e'=T+e$). When the desired conclusion is not immediately reached with the structure of $T'$, we weave together $Z(T')$ and $\dim(T')$ with those of $T$: Bear in mind that $T$ and $T'$ are subject to the same inequalities in passing from being (distinct) trees to the same unicyclic graph, and bear in mind that $\sigma(T')$ may be different from $\sigma(T)$ while $ex(T')$ and $ex(T)$ equal; such a juxtaposition will then yield our desired conclusion. Let $C$ denote the unique cycle in $T+e$. We consider four cases.\\

\emph{Case 1. $C$ contains no emv of $T$:} Let $v$ be an emv of $T$, and let $\{\ell_1, \ell_2, \ldots, \ell_j\}$ be the set of terminal vertices of $v$ in $T$. Further, let $s$ and $s'$ be degree two vertices lying on the path from $\ell_1$ to $v$ in $T$. There are two apparently distinct cases: (A) $e=s\ell_1$ (see (a) of Figure~\ref{fig4}) and (B) $e=ss'$ (see (b) of Figure~\ref{fig4}). In either case, it is immediately clear that $P(T+e)=P(T)+1$, and hence $Z(T+e) \ge Z(T)+1$ by Theorem~\ref{pathcover}. Thus, $\dim(T+e) \le Z(T+e)$ by Theorem~\ref{unicyclic}.

\begin{figure}[ht]
\centering
\begin{tikzpicture}[scale=.65, transform shape]
\node [draw, shape=circle] (a0) at  (-0.5,2) {};
\node [draw, shape=circle] (a1) at  (-1.8,1) {};
\node [draw, shape=circle] (a2) at  (-1.8,0) {};
\node [draw, shape=circle] (a3) at  (-1.8,-1) {};
\node [draw, shape=circle] (a4) at  (-1.8,-2.3) {};
\node [draw, shape=circle] (a11) at  (-0.5,1) {};
\node [draw, shape=circle] (a22) at  (-0.5,0) {};
\node [draw, shape=circle] (a111) at  (0.8,1) {};

\node [draw, shape=circle] (aa0) at  (3.5,2) {};
\node [draw, shape=circle] (aa1) at  (2.2,1) {};
\node [draw, shape=circle] (aa2) at  (2.2,0) {};
\node [draw, shape=circle] (aa11) at  (3.5,1) {};
\node [draw, shape=circle] (aa111) at  (4.8,1) {};

\draw(-1,2)--(a0)--(aa0)--(4,2);
\draw[thick, style=dotted](-2,2)--(-1, 2);
\draw[thick, style=dotted](4,2)--(5,2);
\draw(a0)--(a1)--(a2)--(a3);
\draw(a3)--(-1.8, -1.4);
\draw(-1.8,-1.9)--(a4);
\draw[thick, style=dotted](-1.8,-1.4)--(-1.8,-1.9);
\draw(a22)--(a11)--(a0)--(a111);
\draw[thick, style=dotted](-0.15,0.15)--(0.55,0.75);

\draw(aa0)--(aa1)--(aa2);
\draw(aa11)--(aa0)--(aa111);
\draw[thick, style=dotted](3.8,1)--(4.5, 1);
\draw[very thick](a2) .. controls (-1.1,-0.7) and (-1.1, -1.6) .. (a4);

\node [scale=1.4] at (-0.5,2.5) {$v$};
\node [scale=1.4] at (-2.3,0) {$s$};
\node [scale=1.4] at (-2.3,-2.3) {$\ell_1$};
\node [scale=1.4] at (-0.5,-0.5) {$\ell_2$};
\node [scale=1.4] at (0.9,0.5) {$\ell_j$};
\node [scale=1.4] at (1.5,-3) {\large (a) $T+e$};

\node [draw, shape=circle] (b0) at  (9.5,2) {};
\node [draw, shape=circle] (b1) at  (8.2,1) {};
\node [draw, shape=circle] (b2) at  (8.2,0) {};
\node [draw, shape=circle] (b3) at  (8.2,-1.3) {};
\node [draw, shape=circle] (b4) at  (8.2,-2.6) {};
\node [draw, shape=circle] (b11) at  (9.5,1) {};
\node [draw, shape=circle] (b22) at  (9.5,0) {};
\node [draw, shape=circle] (b111) at  (10.8,1) {};

\node [draw, shape=circle] (bb0) at  (13.5,2) {};
\node [draw, shape=circle] (bb1) at  (12.2,1) {};
\node [draw, shape=circle] (bb2) at  (12.2,0) {};
\node [draw, shape=circle] (bb11) at  (13.5,1) {};
\node [draw, shape=circle] (bb111) at  (14.8,1) {};

\draw(9,2)--(b0)--(bb0)--(14,2);
\draw[thick, style=dotted](8,2)--(9, 2);
\draw[thick, style=dotted](14,2)--(15,2);
\draw(b0)--(b1)--(b2);
\draw(b2)--(8.2, -0.4);
\draw(8.2,-0.9)--(b3);
\draw[thick, style=dotted](8.2,-0.4)--(8.2,-0.9);
\draw(b3)--(8.2, -1.7);
\draw(8.2,-2.2)--(b4);
\draw[thick, style=dotted](8.2,-1.7)--(8.2,-2.2);
\draw(b22)--(b11)--(b0)--(b111);
\draw[thick, style=dotted](9.85,0.15)--(10.55,0.75);

\draw(bb0)--(bb1)--(bb2);
\draw(bb11)--(bb0)--(bb111);
\draw[thick, style=dotted](13.8,1)--(14.5, 1);
\draw[very thick](b1) .. controls (8.9,0.6) and (8.9, -0.9) .. (b3);

\node [scale=1.4] at (9.5,2.5) {$v$};
\node [scale=1.4] at (7.7,1) {$s$};
\node [scale=1.4] at (7.7,-1.3) {$s'$};
\node [scale=1.4] at (7.7,-2.6) {$\ell_1$};
\node [scale=1.4] at (9.5,-0.5) {$\ell_2$};
\node [scale=1.4] at (10.9,0.5) {$\ell_j$};
\node [scale=1.4] at (11.5,-3) {\large (b) $T+e$};

\end{tikzpicture}
\caption{Unicyclic graphs $T+e$ such that $C$ contains no exterior major vertex of $T$}\label{fig4}
\end{figure}
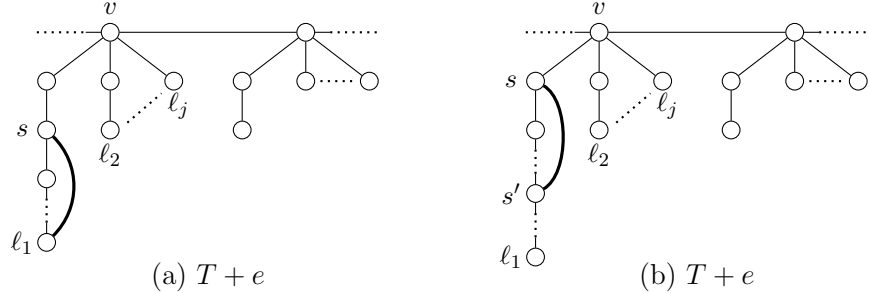

\emph{Case 2. $C$ contains exactly one emv of $T$:} Let $C$ contain one emv, say $v_1$, such that $v_1v_2 \in E(T)$ for some emv $v_2$ with $ter_T(v_2) \ge 2$. Let $\ell_1, \ell_2, \ldots, \ell_j$ be the terminal vertices of $v_1$ and let $s_i$ be a degree two vertex lying on the path from $\ell_i$ to $v_1$ ($1 \le i \le j$). The five a priori cases, as depicted in Figure~\ref{fig5}, reduce to three distinct cases for consideration.

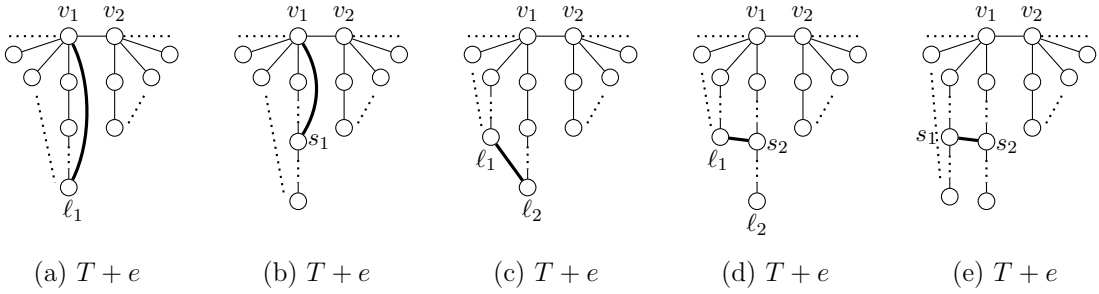
\begin{figure}[ht]
\centering
\begin{tikzpicture}[scale=.61, transform shape]
\node [draw, shape=circle] (a0) at  (0,3) {};
\node [draw, shape=circle] (a1) at  (0,2) {};
\node [draw, shape=circle] (a2) at  (0,1) {};
\node [draw, shape=circle] (a3) at  (0,-0.3) {};
\node [draw, shape=circle] (a11) at  (-0.8,2.1) {};
\node [draw, shape=circle] (a111) at  (-1.2,2.6) {};
\node [draw, shape=circle] (aa0) at  (1,3) {};
\node [draw, shape=circle] (aa1) at  (1,2) {};
\node [draw, shape=circle] (aa2) at  (1,1) {};
\node [draw, shape=circle] (aa11) at  (1.8,2.1) {};
\node [draw, shape=circle] (aa111) at  (2.2,2.6) {};

\draw(a0)--(aa0);
\draw[thick, style=dotted](-1.3,3)--(a0);
\draw[thick, style=dotted](aa0)--(2.3,3);
\draw(a0)--(a1)--(a2);
\draw(a11)--(a0)--(a111);
\draw(a2)--(0, 0.6);
\draw(0,0.1)--(a3);
\draw[thick, style=dotted](0,0.6)--(0,0.1);
\draw(aa0)--(aa1)--(aa2);
\draw(aa11)--(aa0)--(aa111);
\draw[thick, style=dotted](-0.7,1.7)--(-0.3,-0.2);
\draw[thick, style=dotted](1.3,1.1)--(1.7,1.8);
\draw[very thick](a0) .. controls (0.5,2) and (0.5,0.7) .. (a3);

\node [scale=1.4] at (0,3.5) {$v_1$};
\node [scale=1.4] at (1,3.5) {$v_2$};
\node [scale=1.4] at (0.1,-0.8) {$\ell_1$};
\node [scale=1.4] at (0.4,-2.2) {\large (a) $T+e$};

\node [draw, shape=circle] (b0) at  (5,3) {};
\node [draw, shape=circle] (b1) at  (5,2) {};
\node [draw, shape=circle] (b2) at  (5,0.7) {};
\node [draw, shape=circle] (b3) at  (5,-0.6) {};
\node [draw, shape=circle] (b11) at  (4.2,2.1) {};
\node [draw, shape=circle] (b111) at  (3.8,2.6) {};
\node [draw, shape=circle] (bb0) at  (6,3) {};
\node [draw, shape=circle] (bb1) at  (6,2) {};
\node [draw, shape=circle] (bb2) at  (6,1) {};
\node [draw, shape=circle] (bb11) at  (6.8,2.1) {};
\node [draw, shape=circle] (bb111) at  (7.2,2.6) {};

\draw(b0)--(bb0);
\draw[thick, style=dotted](3.7,3)--(b0);
\draw[thick, style=dotted](bb0)--(7.3,3);
\draw(b0)--(b1);
\draw(b11)--(b0)--(b111);
\draw(b1)--(5,1.6);
\draw(5,1.1)--(b2);
\draw[thick, style=dotted](5,1.1)--(5,1.6);
\draw(b2)--(5,0.3);
\draw(5,-0.2)--(b3);
\draw[thick, style=dotted](5,0.3)--(5,-0.2);

\draw(bb0)--(bb1)--(bb2);
\draw(bb11)--(bb0)--(bb111);
\draw[thick, style=dotted](4.3,1.7)--(4.7,-0.5);
\draw[thick, style=dotted](6.3,1.1)--(6.7,1.8);
\draw[very thick](b0) .. controls (5.5,2.2) and (5.5,1.5) .. (b2);

\node [scale=1.4] at (5,3.5) {$v_1$};
\node [scale=1.4] at (6,3.5) {$v_2$};
\node [scale=1.4] at (5.45,0.7) {$s_1$};
\node [scale=1.4] at (5.4,-2.2) {\large (b) $T+e$};

\node [draw, shape=circle] (c0) at  (10,3) {};
\node [draw, shape=circle] (c1) at  (10,2) {};
\node [draw, shape=circle] (c2) at  (10,1) {};
\node [draw, shape=circle] (c3) at  (10,-0.3) {};
\node [draw, shape=circle] (c11) at  (9.2,2.1) {};
\node [draw, shape=circle] (c22) at  (9.2,0.8) {};
\node [draw, shape=circle] (c111) at  (8.8,2.6) {};
\node [draw, shape=circle] (cc0) at  (11,3) {};
\node [draw, shape=circle] (cc1) at  (11,2) {};
\node [draw, shape=circle] (cc2) at  (11,1) {};
\node [draw, shape=circle] (cc11) at  (11.8,2.1) {};
\node [draw, shape=circle] (cc111) at  (12.2,2.6) {};

\draw(c0)--(cc0);
\draw[thick, style=dotted](8.7,3)--(c0);
\draw[thick, style=dotted](cc0)--(12.3,3);
\draw(c0)--(c1)--(c2);
\draw(c11)--(c0)--(c111);
\draw(c2)--(10,0.6);
\draw(10,0.1)--(c3);
\draw[thick, style=dotted](10,0.6)--(10,0.1);
\draw(c11)--(9.2,1.7);
\draw(9.2,1.2)--(c22);
\draw[thick, style=dotted](9.2,1.2)--(9.2,1.7);

\draw(cc0)--(cc1)--(cc2);
\draw(cc11)--(cc0)--(cc111);
\draw[thick, style=dotted](8.8,2.2)--(9,1);
\draw[thick, style=dotted](11.3,1.1)--(11.7,1.8);
\draw[very thick](c3)--(c22);

\node [scale=1.4] at (10,3.5) {$v_1$};
\node [scale=1.4] at (11,3.5) {$v_2$};
\node [scale=1.4] at (10.1,-0.8) {$\ell_2$};
\node [scale=1.4] at (9.1,0.3) {$\ell_1$};
\node [scale=1.4] at (10.4,-2.2) {\large (c) $T+e$};

\node [draw, shape=circle] (d0) at  (15,3) {};
\node [draw, shape=circle] (d1) at  (15,2) {};
\node [draw, shape=circle] (d2) at  (15,0.7) {};
\node [draw, shape=circle] (d3) at  (15,-0.6) {};
\node [draw, shape=circle] (d11) at  (14.2,2.1) {};
\node [draw, shape=circle] (d22) at  (14.2,0.8) {};
\node [draw, shape=circle] (d111) at  (13.8,2.6) {};
\node [draw, shape=circle] (dd0) at  (16,3) {};
\node [draw, shape=circle] (dd1) at  (16,2) {};
\node [draw, shape=circle] (dd2) at  (16,1) {};
\node [draw, shape=circle] (dd11) at  (16.8,2.1) {};
\node [draw, shape=circle] (dd111) at  (17.2,2.6) {};

\draw(d0)--(dd0);
\draw[thick, style=dotted](13.7,3)--(d0);
\draw[thick, style=dotted](dd0)--(17.3,3);
\draw(d0)--(d1);
\draw(d1)--(15,1.6);
\draw(15,1.1)--(d2);
\draw[thick, style=dotted](15,1.6)--(15,1.1);
\draw(d2)--(15,0.3);
\draw(15,-0.2)--(d3);
\draw[thick, style=dotted](15,0.3)--(15,-0.2);

\draw(d11)--(d0)--(d111);
\draw(d11)--(14.2,1.7);
\draw(14.2,1.2)--(d22);
\draw[thick, style=dotted](14.2,1.7)--(14.2,1.2);
\draw[very thick](d2)--(d22);

\draw(dd0)--(dd1)--(dd2);
\draw(dd11)--(dd0)--(dd111);
\draw[thick, style=dotted](13.8,2.2)--(14,1);
\draw[thick, style=dotted](16.3,1.1)--(16.7,1.8);

\node [scale=1.4] at (15,3.5) {$v_1$};
\node [scale=1.4] at (16,3.5) {$v_2$};
\node [scale=1.4] at (15.45,0.6) {$s_2$};
\node [scale=1.4] at (15,-1.1) {$\ell_2$};
\node [scale=1.4] at (14.15,0.3) {$\ell_1$};
\node [scale=1.4] at (15.4,-2.2) {\large (d) $T+e$};

\node [draw, shape=circle] (e0) at  (20,3) {};
\node [draw, shape=circle] (e1) at  (20,2) {};
\node [draw, shape=circle] (e2) at  (20,0.7) {};
\node [draw, shape=circle] (e3) at  (20,-0.6) {};
\node [draw, shape=circle] (e11) at  (19.2,2.1) {};
\node [draw, shape=circle] (e22) at  (19.2,0.8) {};
\node [draw, shape=circle] (e33) at  (19.2,-0.5) {};
\node [draw, shape=circle] (e111) at  (18.8,2.6) {};
\node [draw, shape=circle] (ee0) at  (21,3) {};
\node [draw, shape=circle] (ee1) at  (21,2) {};
\node [draw, shape=circle] (ee2) at  (21,1) {};
\node [draw, shape=circle] (ee11) at  (21.8,2.1) {};
\node [draw, shape=circle] (ee111) at  (22.2,2.6) {};

\draw(e0)--(ee0);
\draw[thick, style=dotted](18.7,3)--(e0);
\draw[thick, style=dotted](ee0)--(22.3,3);
\draw(e0)--(e1);
\draw(e1)--(20,1.6);
\draw(20,1.1)--(e2);
\draw[thick, style=dotted](20,1.6)--(20,1.1);
\draw(e2)--(20,0.3);
\draw(20,-0.2)--(e3);
\draw[thick, style=dotted](20,0.3)--(20,-0.2);

\draw(e11)--(e0)--(e111);
\draw(e11)--(19.2,1.7);
\draw(19.2,1.2)--(e22);
\draw[thick, style=dotted](19.2,1.7)--(19.2,1.2);
\draw(e22)--(19.2,0.4);
\draw(19.2,-0.1)--(e33);
\draw[thick, style=dotted](19.2,0.4)--(19.2,-0.1);

\draw[very thick](e2)--(e22);

\draw(ee0)--(ee1)--(ee2);
\draw(ee11)--(ee0)--(ee111);
\draw[thick, style=dotted](18.8,2.2)--(19,-0.2);
\draw[thick, style=dotted](21.3,1.1)--(21.7,1.8);

\node [scale=1.4] at (20,3.5) {$v_1$};
\node [scale=1.4] at (21,3.5) {$v_2$};
\node [scale=1.4] at (20.45,0.6) {$s_2$};
\node [scale=1.4] at (18.7,0.8) {$s_1$};
\node [scale=1.4] at (20.4,-2.2) {\large (e) $T+e$};

\end{tikzpicture}
\caption{Unicyclic graphs $T+e$ such that $C$ contains exactly one exterior major vertex of $T$}\label{fig5}
\end{figure}

\emph{Subcase 2.1. $e=v_1\ell_1$ or $e=\ell_1\ell_2$ ((a) or (c), respectively, of Figure~\ref{fig5}):} By removing an edge from $T+e$ in (a) of Figure~\ref{fig5}, one obtains a tree $T$ in (c) of Figure~\ref{fig5}. So, we only need to consider the case $e=\ell_1\ell_2$. Let $x \in N(v_1) \cap V(C)$ ($C$ here is the unique cycle), and let $e'=v_1x$. If $Z(T')>\dim(T')$, then we're done; otherwise we have $Z(T')=\dim(T')=\sigma(T')-ex(T')$. Suppose $ter_T(v_1) \ge 3$. Then $ex(T')=ex(T)$, whereas $\sigma(T')=\sigma(T)-1$ since $ter_{T'}(v_1)=ter_{T}(v_1)-1$; thus $Z(T')=Z(T)-1$. And we have $Z(T'+e')=Z(T+e) \ge Z(T)-1=Z(T')$. If $ter_T(v_1)=2$ and $v_1$ is adjacent to at least two emvs in $T$, then $v_1$ becomes a major vertex with $ter_{T'}(v_1)=1$. If $ter_T(v_1)=2$ and $v_1$ is adjacent to exactly one emv $v_2$ in $T$, then $v_1$ becomes an exterior degree two vertex in $T'$, and we have $ex(T')=k$.

\emph{Subcase 2.2. $e=v_1s_1$ or $e=\ell_1s_2$ ((b) or (d), respectively, of Figure~\ref{fig5}):} Notice (b) is a special case of (d) when $d_T(v_1, s_2)=1$. So, we only need to consider the case $e=\ell_1s_2$. Let $x$ be the vertex adjacent to $v_1$ and lying on the path from $\ell_2$ to $v_1$ in $T$. Take $e'=v_1x$. If $d_T(v_1, s_2) \ge 2$, then $\ell_1$ becomes an interior degree two vertex in $T'$. Next, suppose $d_T(v_1, s_2)=1$. If $Z(T')>\dim(T')$, then we're done; otherwise we have $Z(T')=\dim(T')=\sigma(T')-ex(T')$. If $ter_T(v_1) \ge 3$, then $ex(T')=ex(T)$, whereas $\sigma(T')=\sigma(T)-1$ since $ter_{T'}(v_1)=ter_T(v_1)-1$; thus $Z(T')=Z(T)-1$. And we have $Z(T'+e')=Z(T+e) \ge Z(T)-1=Z(T')$. If $ter_T(v_1)=2$ and $v_1$ is adjacent to at least two emvs in $T$, then $v_1$ becomes an emv with $ter_{T'}(v_1)=1$. If $ter_T(v_1)=2$ and $v_1$ is adjacent to exactly one emv $v_2$ in $T$, then $v_1$ becomes an exterior degree two vertex in $T'$, and thus $ex(T')=k$.

\emph{Subcase 2.3. $e=s_1s_2$ (see (e) of Figure~\ref{fig5}):} Suppose $d_T(v_1, s_1) \ge 2$ or $d_T(v_1, s_2) \ge 2$; assume WLOG the former, take $e'=v_1x$, where $x \in N(v_1)$ lies on the path from $\ell_1$ to $v_1$ in $T$. Then $s_2$ becomes an emv in $T'$ with $ter_{T'}(s_2)=1$. Next, suppose that $d_T(v_1, s_1)=1=d_T(v_1, s_2)$, and we take $e'=v_1s_1$. If $ter_T(v_1) \ge 4$, then $ex(T')=ex(T)+1$ and $\sigma(T')=\sigma(T)$, and hence $\dim(T')=\dim(T)-1$; thus $\dim(T+e)=\dim(T'+e') \le \dim(T')+1=\dim(T)$. If $ter_T(v_1)=3$, then $v_1$ becomes an emv with $ter_{T'}(v_1)=1$. If $ter_T(v_1)=2$ and $v_1$ is adjacent to at least two emvs in $T$, then $v_1$ becomes a major vertex with $ter_{T'}(v_1)=0$. If $ter_T(v_1)=2$ and $v_1$ is adjacent to exactly one emv $v_2$ in $T$, then $v_1$ becomes an interior degree two vertex in $T'$.\\

\emph{Case 3. $C$ contains exactly two emvs of $T$:} Let $C$ contain two emvs, say $v_1, v_2$, such that $v_1v_2 \in E(C)$. For each $v_i$ ($i=1,2$), let $\ell_i$ be a terminal vertex of $v_i$ and let $s_i$ be a degree two vertex lying on the path from $\ell_i$ to $v_i$. We consider five subcases.

\begin{figure}[ht]
\centering
\begin{tikzpicture}[scale=.61, transform shape]
\node [draw, shape=circle] (a0) at  (0,3) {};
\node [draw, shape=circle] (a1) at  (0,2) {};
\node [draw, shape=circle] (a2) at  (0,0.7) {};
\node [draw, shape=circle] (a11) at  (-0.8,2.1) {};
\node [draw, shape=circle] (a111) at  (-1.2,2.6) {};
\node [draw, shape=circle] (aa0) at  (1,3) {};
\node [draw, shape=circle] (aa1) at  (1,2) {};
\node [draw, shape=circle] (aa2) at  (1,1) {};
\node [draw, shape=circle] (aa3) at  (1,-0.3) {};
\node [draw, shape=circle] (aa11) at  (1.8,2.1) {};
\node [draw, shape=circle] (aa111) at  (2.2,2.6) {};

\draw(a0)--(aa0);
\draw[thick, style=dotted](-1.3,3)--(a0);
\draw[thick, style=dotted](aa0)--(2.3,3);
\draw(a0)--(a1);
\draw(a1)--(0,1.6);
\draw(0,1.1)--(a2);
\draw[thick, style=dotted](0,1.6)--(0,1.1);

\draw(a11)--(a0)--(a111);
\draw(aa0)--(aa1)--(aa2);
\draw(aa2)--(1,0.6);
\draw(1,0.1)--(aa3);
\draw[thick, style=dotted](1,0.6)--(1,0.1);
\draw(aa11)--(aa0)--(aa111);
\draw[thick, style=dotted](-0.7,1.7)--(-0.3,0.7);
\draw[thick, style=dotted](1.3,0)--(1.7,1.7);
\draw[very thick](a0)--(aa3);

\node [scale=1.4] at (0,3.5) {$v_1$};
\node [scale=1.4] at (1,3.5) {$v_2$};
\node [scale=1.4] at (1,-0.8) {$\ell_2$};
\node [scale=1.4] at (0.4,-2.2) {\large (a) $T+e$};

\node [draw, shape=circle] (b0) at  (5,3) {};
\node [draw, shape=circle] (b1) at  (5,2) {};
\node [draw, shape=circle] (b2) at  (5,0.7) {};
\node [draw, shape=circle] (b11) at  (4.2,2.1) {};
\node [draw, shape=circle] (b111) at  (3.8,2.6) {};
\node [draw, shape=circle] (bb0) at  (6,3) {};
\node [draw, shape=circle] (bb1) at  (6,2) {};
\node [draw, shape=circle] (bb2) at  (6,0.7) {};
\node [draw, shape=circle] (bb3) at  (6,-0.6) {};
\node [draw, shape=circle] (bb11) at  (6.8,2.1) {};
\node [draw, shape=circle] (bb111) at  (7.2,2.6) {};

\draw(b0)--(bb0);
\draw[thick, style=dotted](3.7,3)--(b0);
\draw[thick, style=dotted](bb0)--(7.3,3);
\draw(b0)--(b1);
\draw(b11)--(b0)--(b111);
\draw(b1)--(5,1.6);
\draw(5,1.1)--(b2);
\draw[thick, style=dotted](5,1.1)--(5,1.6);

\draw(bb0)--(bb1);
\draw(bb11)--(bb0)--(bb111);
\draw(bb1)--(6,1.6);
\draw(6,1.1)--(bb2);
\draw[thick, style=dotted](6,1.6)--(6,1.1);
\draw(bb2)--(6,0.3);
\draw(6,-0.2)--(bb3);
\draw[thick, style=dotted](6,0.3)--(6,-0.2);
\draw[thick, style=dotted](4.3,1.7)--(4.7,0.7);
\draw[thick, style=dotted](6.3,-0.3)--(6.7,1.7);
\draw[very thick](b0)--(bb2);

\node [scale=1.4] at (5,3.5) {$v_1$};
\node [scale=1.4] at (6,3.5) {$v_2$};
\node [scale=1.4] at (5.55,0.7) {$s_2$};
\node [scale=1.4] at (5.4,-2.2) {\large (b) $T+e$};

\node [draw, shape=circle] (c0) at  (10,3) {};
\node [draw, shape=circle] (c1) at  (10,2) {};
\node [draw, shape=circle] (c2) at  (10,0.7) {};
\node [draw, shape=circle] (c11) at  (9.2,2.1) {};
\node [draw, shape=circle] (c111) at  (8.8,2.6) {};
\node [draw, shape=circle] (cc0) at  (11,3) {};
\node [draw, shape=circle] (cc1) at  (11,2) {};
\node [draw, shape=circle] (cc2) at  (11,1) {};
\node [draw, shape=circle] (cc3) at  (11,-0.3) {};
\node [draw, shape=circle] (cc11) at  (11.8,2.1) {};
\node [draw, shape=circle] (cc111) at  (12.2,2.6) {};

\draw(c0)--(cc0);
\draw[thick, style=dotted](8.7,3)--(c0);
\draw[thick, style=dotted](cc0)--(12.3,3);
\draw(c11)--(c0)--(c111);
\draw(c0)--(c1);
\draw(c1)--(10,1.6);
\draw(10,1.1)--(c2);
\draw[thick, style=dotted](10,1.1)--(10,1.6);
\draw(cc0)--(cc1)--(cc2);
\draw(cc11)--(cc0)--(cc111);
\draw(cc2)--(11,0.6);
\draw(11,0.1)--(cc3);
\draw[thick, style=dotted](11,0.6)--(11,0.1);
\draw[thick, style=dotted](9.3,1.7)--(9.7,0.8);
\draw[thick, style=dotted](11.3,-0.1)--(11.7,1.7);
\draw[very thick](c2)--(cc3);

\node [scale=1.4] at (10,3.5) {$v_1$};
\node [scale=1.4] at (11,3.5) {$v_2$};
\node [scale=1.4] at (10,0.2) {$\ell_1$};
\node [scale=1.4] at (11.1,-0.8) {$\ell_2$};
\node [scale=1.4] at (10.4,-2.2) {\large (c) $T+e$};

\node [draw, shape=circle] (d0) at  (15,3) {};
\node [draw, shape=circle] (d1) at  (15,2) {};
\node [draw, shape=circle] (d2) at  (15,0.7) {};
\node [draw, shape=circle] (d11) at  (14.2,2.1) {};
\node [draw, shape=circle] (d111) at  (13.8,2.6) {};
\node [draw, shape=circle] (dd0) at  (16,3) {};
\node [draw, shape=circle] (dd1) at  (16,2) {};
\node [draw, shape=circle] (dd2) at  (16,0.7) {};
\node [draw, shape=circle] (dd3) at  (16,-0.6) {};
\node [draw, shape=circle] (dd11) at  (16.8,2.1) {};
\node [draw, shape=circle] (dd22) at  (17.1,1.1) {};
\node [draw, shape=circle] (dd111) at  (17.2,2.6) {};

\draw(d0)--(dd0);
\draw[thick, style=dotted](13.7,3)--(d0);
\draw[thick, style=dotted](dd0)--(17.3,3);
\draw(d0)--(d1);
\draw(d11)--(d0)--(d111);
\draw(d1)--(15,1.6);
\draw(15,1.1)--(d2);
\draw[thick, style=dotted](15,1.1)--(15,1.6);

\draw(dd0)--(dd1);
\draw(dd22)--(dd11)--(dd0)--(dd111);
\draw(dd1)--(16,1.6);
\draw(16,1.1)--(dd2);
\draw[thick, style=dotted](16,1.6)--(16,1.1);
\draw(dd2)--(16,0.3);
\draw(16,-0.2)--(dd3);
\draw[thick, style=dotted](16,0.3)--(16,-0.2);
\draw[thick, style=dotted](14.3,1.7)--(14.7,0.7);
\draw[thick, style=dotted](16.3,-0.3)--(17,0.9);
\draw[very thick](d2)--(dd2);

\node [scale=1.4] at (15,3.5) {$v_1$};
\node [scale=1.4] at (16,3.5) {$v_2$};
\node [scale=1.4] at (15,0.2) {$\ell_1$};
\node [scale=1.4] at (16.45,0.7) {$s_2$};
\node [scale=1.4] at (15.4,-2.2) {\large (d) $T+e$};

\node [draw, shape=circle] (e0) at  (20,3) {};
\node [draw, shape=circle] (e1) at  (20,2) {};
\node [draw, shape=circle] (e2) at  (20,0.7) {};
\node [draw, shape=circle] (e3) at  (20,-0.6) {};
\node [draw, shape=circle] (e11) at  (19.2,2.1) {};
\node [draw, shape=circle] (e22) at  (18.9,1.1) {};
\node [draw, shape=circle] (e111) at  (18.8,2.6) {};
\node [draw, shape=circle] (ee0) at  (21,3) {};
\node [draw, shape=circle] (ee1) at  (21,2) {};
\node [draw, shape=circle] (ee2) at  (21,0.7) {};
\node [draw, shape=circle] (ee3) at  (21,-0.6) {};
\node [draw, shape=circle] (ee11) at  (21.8,2.1) {};
\node [draw, shape=circle] (ee22) at  (22.1,1.1) {};
\node [draw, shape=circle] (ee111) at  (22.2,2.6) {};

\draw(e0)--(ee0);
\draw[thick, style=dotted](18.7,3)--(e0);
\draw[thick, style=dotted](ee0)--(22.3,3);
\draw(e0)--(e1);
\draw(e22)--(e11)--(e0)--(e111);
\draw(e1)--(20,1.6);
\draw(20,1.1)--(e2);
\draw[thick, style=dotted](20,1.1)--(20,1.6);
\draw(e2)--(20,0.3);
\draw(20,-0.2)--(e3);
\draw[thick, style=dotted](20,0.3)--(20,-0.2);

\draw(ee0)--(ee1);
\draw(ee22)--(ee11)--(ee0)--(ee111);
\draw(ee1)--(21,1.6);
\draw(21,1.1)--(ee2);
\draw[thick, style=dotted](21,1.6)--(21,1.1);
\draw(ee2)--(21,0.3);
\draw(21,-0.2)--(ee3);
\draw[thick, style=dotted](21,0.3)--(21,-0.2);
\draw[thick, style=dotted](19,0.8)--(19.7,-0.3);
\draw[thick, style=dotted](21.3,-0.3)--(22,0.8);
\draw[very thick](e2)--(ee2);

\node [scale=1.4] at (20,3.5) {$v_1$};
\node [scale=1.4] at (21,3.5) {$v_2$};
\node [scale=1.4] at (19.55,0.7) {$s_1$};
\node [scale=1.4] at (21.45,0.7) {$s_2$};
\node [scale=1.4] at (20.4,-2.2) {\large (e) $T+e$};

\end{tikzpicture}
\caption{Unicyclic graphs $T+e$ such that $C$ contains exactly two exterior major vertices of $T$}\label{fig6}
\end{figure}
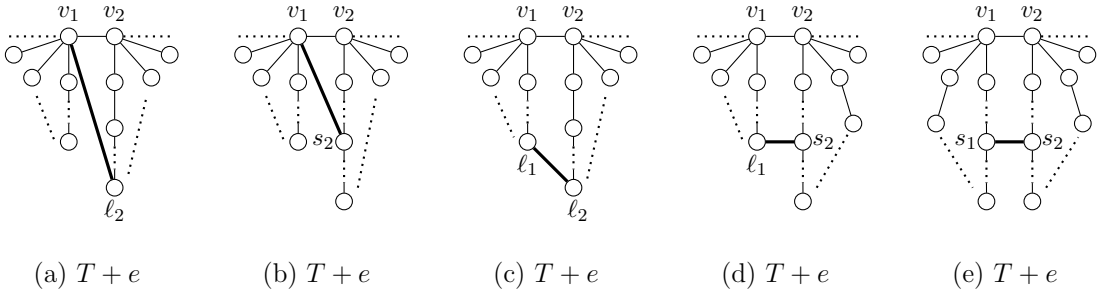

\emph{Subcase 3.1. $e=v_1\ell_2$ (see (a) of Figure~\ref{fig6}):} Let $e'=v_1v_2$. If $ter_T(v_2) \ge 3$ or $v_2$ is adjacent to at least two emvs in $T$, then $\ell_2$ becomes an interior degree two vertex in $T'$. If $ter_T(v_2)=2$ and $v_2$ is adjacent to exactly one emv in $T$, then $v_2$ becomes an exterior degree two vertex in $T'$ without turning any non-emv of $T$ into an emv of $T'$, and thus $ex(T')=k$.

\emph{Subcase 3.2. $e=v_1s_2$ (see (b) of Figure~\ref{fig6}):} If $ter_T(v_2) \ge 3$ or $v_2$ is adjacent to at least two emvs in $T$, then take $e'=v_1v_2$; notice that $s_2$ becomes an emv of terminal degree one in $T'$. Next, we consider when $ter_T(v_2)=2$ and $v_2$ is adjacent to exactly one emv in $T$. Let $u$ be the vertex adjacent to $s_2$ and lying on the $s_2-v_2$ path of $T$. Now, take $e'=s_2u$. If $u=v_2$, then $ex(T')=k$. If $u\neq v_2$, then $\sigma(T')=\sigma(T)+1$, while $ex(T')=ex(T)$. Again, if $Z(T')>\dim(T')$, then we're done. Otherwise, $Z(T')=\dim(T')=\dim(T)+1$. And $Z(T+e)=Z(T'+e')\geq Z(T')-1$ implies that $Z(T+e)+1\geq\dim(T)+1\geq\dim(T+e)$.

\emph{Subcase 3.3. $e=\ell_1\ell_2$ (see (c) of Figure~\ref{fig6}):} Let $e'=v_1v_2$. If $ter_T(v_i) \ge 3$ or $v_i$ is adjacent to at least two emvs in $T$ (i.e., $\deg_T(v_i) \ge 4$) for each $i=1,2$, then $\ell_1$ and $\ell_2$ become interior degree two vertices in $T'$. If $\deg_T(v_i) \ge 4$ for either $i=1$ or $i=2$ but not both, say $\deg_T(v_1) \ge 4$ and $\deg_T(v_2)=3$, then $v_2$ becomes an exterior degree two vertex in $T'$ without increasing the number of emvs, and hence $ex(T')=k$. If $\deg_T(v_1)=\deg_T(v_2)=3$, then $T'$ becomes a path, and thus, $\dim(T'+e')=Z(T'+e')=2$.

\emph{Subcase 3.4. $e=\ell_1s_2$ (see (d) of Figure~\ref{fig6}):} Take $e'=v_1v_2$. The only way for $\ell_1$ to not be an interior degree $2$ vertex in $T'$ is for $\deg_{T}(v_1)=3$; thus, $ter_{T'}(s_2)\geq 2$. If $ter_{T'}(v_2)=0$, then $ex(T')=k$. If $ter_{T'}(v_2)=1$, then $\dim(T') < Z(T')$ by Theorem \ref{dim=Z,T}. If $ter_{T'}(v_2)\geq 2$, then we have $\sigma(T')=\sigma(T)-1$, while $ex(T')=ex(T)$. And the same argument (transposing $T$ and $T'$) as in subcase 3.2 applies.

\emph{Subcase 3.5. $e=s_1s_2$ (see (e) of Figure~\ref{fig6}):} Take $e'=v_1v_2$; then $\deg_{T'}(s_i)=3$ for $i=1,2$. The only way for $ter_{T'}(s_i)\geq 2$ ($i=1,2$) is for $\deg_{T}(v_1)=3=\deg_{T}(v_2)$. For such a tree $T$, it's obvious that $Z(T+e)\geq 2=Z(T)$.\\

\emph{Case 4. $C$ contains three or more emvs of $T$:} Let $C$ contain $t \ge 3$ emvs, say $v_1, v_2, \ldots, v_t$, such that $v_iv_{i+1} \in E(C)$, where $1 \le i \le t-1$. For each $v_i$ ($1 \le i \le t$), let $\ell_i$ be a terminal vertex of $v_i$ and let $s_i$ be a degree two vertex lying on the path from $\ell_i$ to $v_i$. We consider six subcases.

\emph{Subcase 4.1. $C$ contains only emvs of $T$:} Let a set $W$ contain all but one of the terminal vertices of $v_i$ for every emv $v_i$ of $T$; it is readily checked that $W$ forms a resolving set for $T+e$ (c.f. Case 2 of Proposition~\ref{ex1}), and thus $\dim(T+e) \le \dim(T)$.

\begin{figure}[ht]
\centering
\begin{tikzpicture}[scale=.495, transform shape]

\node [draw, shape=circle] (a0) at  (0,3) {};
\node [draw, shape=circle] (a1) at  (0,2) {};
\node [draw, shape=circle] (a2) at  (0,0.7) {};
\node [draw, shape=circle] (a11) at  (-0.8,2.1) {};
\node [draw, shape=circle] (a111) at  (-1.2,2.6) {};
\node [draw, shape=circle] (aa0) at  (1,3) {};
\node [draw, shape=circle] (aa1) at  (0.5,3.9) {};
\node [draw, shape=circle] (aa11) at  (1,4.3) {};
\node [draw, shape=circle] (aa111) at  (1.8,3.8) {};
\node [draw, shape=circle] (aaa0) at  (2.5,3) {};
\node [draw, shape=circle] (aaa1) at  (2.5,2) {};
\node [draw, shape=circle] (aaa2) at  (2.5,1) {};
\node [draw, shape=circle] (aaa3) at  (2.5,-0.3) {};
\node [draw, shape=circle] (aaa11) at  (3.3,2.1) {};
\node [draw, shape=circle] (aaa111) at  (3.7,2.6) {};

\draw[thick, style=dotted](-1.3,3)--(a0);
\draw(a0)--(aa0);
\draw(aa0)--(1.4,3);
\draw(2.1,3)--(aaa0);
\draw[thick, style=dotted](1.4,3)--(2.1,3);
\draw[thick, style=dotted](aaa0)--(3.8,3);

\draw(a0)--(a1);
\draw(a1)--(0,1.6);
\draw(0,1.1)--(a2);
\draw[thick, style=dotted](0,1.6)--(0,1.1);

\draw(a11)--(a0)--(a111);
\draw(aa1)--(aa0)--(aa11);
\draw(aa0)--(aa111);
\draw(aaa0)--(aaa1)--(aaa2);
\draw(aaa2)--(2.5,0.6);
\draw(2.5,0.1)--(aaa3);
\draw[thick, style=dotted](2.5,0.6)--(2.5,0.1);
\draw(aaa11)--(aaa0)--(aaa111);
\draw[thick, style=dotted](-0.7,1.7)--(-0.3,0.7);
\draw[thick, style=dotted](2.8,0)--(3.2,1.7);
\draw[thick, style=dotted](1.25,4.2)--(1.6,3.95);
\draw[very thick](a0)--(aaa3);

\node [scale=1.4] at (0,3.45) {$v_1$};
\node [scale=1.4] at (1,2.55) {$v_2$};
\node [scale=1.4] at (2.5,3.45) {$v_t$};
\node [scale=1.4] at (2.5,-0.8) {$\ell_t$};
\node [scale=1.4] at (1.1,-2.2) {\large (a) $T+e$};

\node [draw, shape=circle] (b0) at  (6,3) {};
\node [draw, shape=circle] (b1) at  (6,2) {};
\node [draw, shape=circle] (b2) at  (6,0.7) {};
\node [draw, shape=circle] (b11) at  (5.2,2.1) {};
\node [draw, shape=circle] (b111) at  (4.8,2.6) {};
\node [draw, shape=circle] (bb0) at  (7,3) {};
\node [draw, shape=circle] (bb1) at  (6.5,3.9) {};
\node [draw, shape=circle] (bb11) at  (7,4.3) {};
\node [draw, shape=circle] (bb111) at  (7.8,3.8) {};
\node [draw, shape=circle] (bbb0) at  (8.5,3) {};
\node [draw, shape=circle] (bbb1) at  (8.5,2) {};
\node [draw, shape=circle] (bbb2) at  (8.5,0.7) {};
\node [draw, shape=circle] (bbb3) at  (8.5,-0.6) {};
\node [draw, shape=circle] (bbb11) at  (9.3,2.1) {};
\node [draw, shape=circle] (bbb111) at  (9.7,2.6) {};

\draw[thick, style=dotted](4.7,3)--(b0);
\draw(b0)--(bb0);
\draw(bb0)--(7.4,3);
\draw(8.1,3)--(bbb0);
\draw[thick, style=dotted](7.4,3)--(8.1,3);
\draw[thick, style=dotted](bbb0)--(9.8,3);

\draw(b0)--(b1);
\draw(b1)--(6,1.6);
\draw(6,1.1)--(b2);
\draw[thick, style=dotted](6,1.6)--(6,1.1);

\draw(b11)--(b0)--(b111);
\draw(bb1)--(bb0)--(bb11);
\draw(bb0)--(bb111);
\draw(bbb0)--(bbb1);
\draw(bbb1)--(8.5,1.6);
\draw(8.5,1.1)--(bbb2);
\draw[thick, style=dotted](8.5,1.6)--(8.5,1.1);
\draw(bbb2)--(8.5,0.3);
\draw(8.5,-0.2)--(bbb3);
\draw[thick, style=dotted](8.5,0.3)--(8.5,-0.2);
\draw(bbb11)--(bbb0)--(bbb111);
\draw[thick, style=dotted](5.3,1.7)--(5.7,0.7);
\draw[thick, style=dotted](8.8,-0.3)--(9.2,1.7);
\draw[thick, style=dotted](7.25,4.2)--(7.6,3.95);
\draw[very thick](b0)--(bbb2);

\node [scale=1.4] at (6,3.45) {$v_1$};
\node [scale=1.4] at (7,2.55) {$v_2$};
\node [scale=1.4] at (8.5,3.45) {$v_t$};
\node [scale=1.4] at (8.05,0.7) {$s_t$};
\node [scale=1.4] at (7.1,-2.2) {\large (b) $T+e$};

\node [draw, shape=circle] (c0) at  (12,3) {};
\node [draw, shape=circle] (c1) at  (12,2) {};
\node [draw, shape=circle] (c2) at  (12,0.7) {};
\node [draw, shape=circle] (c11) at  (11.2,2.1) {};
\node [draw, shape=circle] (c111) at  (10.8,2.6) {};
\node [draw, shape=circle] (cc0) at  (13,3) {};
\node [draw, shape=circle] (cc1) at  (12.5,3.9) {};
\node [draw, shape=circle] (cc11) at  (13,4.3) {};
\node [draw, shape=circle] (cc111) at  (13.8,3.8) {};
\node [draw, shape=circle] (ccc0) at  (14.5,3) {};
\node [draw, shape=circle] (ccc1) at  (14.5,2) {};
\node [draw, shape=circle] (ccc2) at  (14.5,1) {};
\node [draw, shape=circle] (ccc3) at  (14.5,-0.3) {};
\node [draw, shape=circle] (ccc11) at  (15.3,2.1) {};
\node [draw, shape=circle] (ccc111) at  (15.7,2.6) {};

\draw[thick, style=dotted](10.7,3)--(c0);
\draw(c0)--(cc0);
\draw(cc0)--(13.4,3);
\draw(14.1,3)--(ccc0);
\draw[thick, style=dotted](13.4,3)--(14.1,3);
\draw[thick, style=dotted](ccc0)--(15.8,3);

\draw(c0)--(c1);
\draw(c1)--(12,1.6);
\draw(12,1.1)--(c2);
\draw[thick, style=dotted](12,1.6)--(12,1.1);

\draw(c11)--(c0)--(c111);
\draw(cc1)--(cc0)--(cc11);
\draw(cc0)--(cc111);
\draw(ccc0)--(ccc1)--(ccc2);
\draw(ccc2)--(14.5,0.6);
\draw(14.5,0.1)--(ccc3);
\draw[thick, style=dotted](14.5,0.6)--(14.5,0.1);
\draw(ccc11)--(ccc0)--(ccc111);
\draw[thick, style=dotted](11.3,1.7)--(11.7,0.7);
\draw[thick, style=dotted](14.8,0)--(15.2,1.7);
\draw[thick, style=dotted](13.25,4.2)--(13.6,3.95);
\draw[very thick](c2)--(ccc3);

\node [scale=1.4] at (12,3.45) {$v_1$};
\node [scale=1.4] at (13,2.55) {$v_2$};
\node [scale=1.4] at (14.5,3.45) {$v_t$};
\node [scale=1.4] at (12.05,0.2) {$\ell_1$};
\node [scale=1.4] at (14.55,-0.8) {$\ell_t$};
\node [scale=1.4] at (13.1,-2.2) {\large (c) $T+e$};

\node [draw, shape=circle] (d0) at  (18,3) {};
\node [draw, shape=circle] (d1) at  (18,2) {};
\node [draw, shape=circle] (d2) at  (18,0.7) {};
\node [draw, shape=circle] (d11) at  (17.2,2.1) {};
\node [draw, shape=circle] (d111) at  (16.8,2.6) {};
\node [draw, shape=circle] (dd0) at  (19,3) {};
\node [draw, shape=circle] (dd1) at  (18.5,3.9) {};
\node [draw, shape=circle] (dd11) at  (19,4.3) {};
\node [draw, shape=circle] (dd111) at  (19.8,3.8) {};
\node [draw, shape=circle] (ddd0) at  (20.5,3) {};
\node [draw, shape=circle] (ddd1) at  (20.5,2) {};
\node [draw, shape=circle] (ddd2) at  (20.5,0.7) {};
\node [draw, shape=circle] (ddd3) at  (20.5,-0.6) {};
\node [draw, shape=circle] (ddd11) at  (21.3,2.1) {};
\node [draw, shape=circle] (ddd111) at  (21.7,2.6) {};

\draw[thick, style=dotted](16.7,3)--(d0);
\draw(d0)--(dd0);
\draw(dd0)--(19.4,3);
\draw(20.1,3)--(ddd0);
\draw[thick, style=dotted](19.4,3)--(20.1,3);
\draw[thick, style=dotted](ddd0)--(21.8,3);

\draw(d0)--(d1);
\draw(d1)--(18,1.6);
\draw(18,1.1)--(d2);
\draw[thick, style=dotted](18,1.6)--(18,1.1);
\draw(d11)--(d0)--(d111);
\draw(dd1)--(dd0)--(dd11);

\draw(dd0)--(dd111);
\draw(ddd0)--(ddd1);
\draw(ddd1)--(20.5,1.6);
\draw(20.5,1.1)--(ddd2);
\draw[thick, style=dotted](20.5,1.6)--(20.5,1.1);
\draw(ddd2)--(20.5,0.3);
\draw(20.5,-0.2)--(ddd3);
\draw[thick, style=dotted](20.5,0.3)--(20.5,-0.2);
\draw(ddd11)--(ddd0)--(ddd111);
\draw[thick, style=dotted](17.3,1.7)--(17.7,0.7);
\draw[thick, style=dotted](20.8,-0.3)--(21.2,1.7);
\draw[thick, style=dotted](19.25,4.2)--(19.6,3.95);
\draw[very thick](d2)--(ddd2);

\node [scale=1.4] at (18,3.45) {$v_1$};
\node [scale=1.4] at (19,2.55) {$v_2$};
\node [scale=1.4] at (20.5,3.45) {$v_t$};
\node [scale=1.4] at (18.05,0.2) {$\ell_1$};
\node [scale=1.4] at (20.2,0.4) {$s_t$};
\node [scale=1.4] at (19.1,-2.2) {\large (d) $T+e$};

\node [draw, shape=circle] (e0) at  (24,3) {};
\node [draw, shape=circle] (e1) at  (24,2) {};
\node [draw, shape=circle] (e2) at  (24,0.7) {};
\node [draw, shape=circle] (e3) at  (24,-0.6) {};
\node [draw, shape=circle] (e11) at  (23.2,2.1) {};
\node [draw, shape=circle] (e111) at  (22.8,2.6) {};
\node [draw, shape=circle] (ee0) at  (25,3) {};
\node [draw, shape=circle] (ee1) at  (24.5,3.9) {};
\node [draw, shape=circle] (ee11) at  (25,4.3) {};
\node [draw, shape=circle] (ee111) at  (25.8,3.8) {};
\node [draw, shape=circle] (eee0) at  (26.5,3) {};
\node [draw, shape=circle] (eee1) at  (26.5,2) {};
\node [draw, shape=circle] (eee2) at  (26.5,0.7) {};
\node [draw, shape=circle] (eee3) at  (26.5,-0.6) {};
\node [draw, shape=circle] (eee11) at  (27.3,2.1) {};
\node [draw, shape=circle] (eee111) at  (27.7,2.6) {};

\draw[thick, style=dotted](22.7,3)--(e0);
\draw(e0)--(ee0);
\draw(ee0)--(25.4,3);
\draw(26.1,3)--(eee0);
\draw[thick, style=dotted](25.4,3)--(26.1,3);
\draw[thick, style=dotted](eee0)--(27.8,3);

\draw(e0)--(e1);
\draw(e1)--(24,1.6);
\draw(24,1.1)--(e2);
\draw[thick, style=dotted](24,1.6)--(24,1.1);
\draw(e2)--(24,0.3);
\draw(24,-0.2)--(e3);
\draw[thick, style=dotted](24,0.3)--(24,-0.2);

\draw(e11)--(e0)--(e111);
\draw(ee1)--(ee0)--(ee11);

\draw(ee0)--(ee111);
\draw(eee0)--(eee1);
\draw(eee1)--(26.5,1.6);
\draw(26.5,1.1)--(eee2);
\draw[thick, style=dotted](26.5,1.6)--(26.5,1.1);
\draw(eee2)--(26.5,0.3);
\draw(26.5,-0.2)--(eee3);
\draw[thick, style=dotted](26.5,0.3)--(26.5,-0.2);
\draw(eee11)--(eee0)--(eee111);
\draw[thick, style=dotted](23.3,1.7)--(23.7,-0.3);%
\draw[thick, style=dotted](26.8,-0.3)--(27.2,1.7);
\draw[thick, style=dotted](25.25,4.2)--(25.6,3.95);
\draw[very thick](e2)--(eee2);

\node [scale=1.4] at (24,3.45) {$v_1$};
\node [scale=1.4] at (25,2.55) {$v_2$};
\node [scale=1.4] at (26.5,3.45) {$v_t$};
\node [scale=1.4] at (24.4,0.4) {$s_1$};
\node [scale=1.4] at (26.2,0.4) {$s_t$};
\node [scale=1.4] at (25.1,-2.2) {\large (e) $T+e$};

\end{tikzpicture}
\caption{Unicyclic graphs $T+e$ such that $C$ contains at least three exterior major vertices of $T$}\label{fig7}
\end{figure}
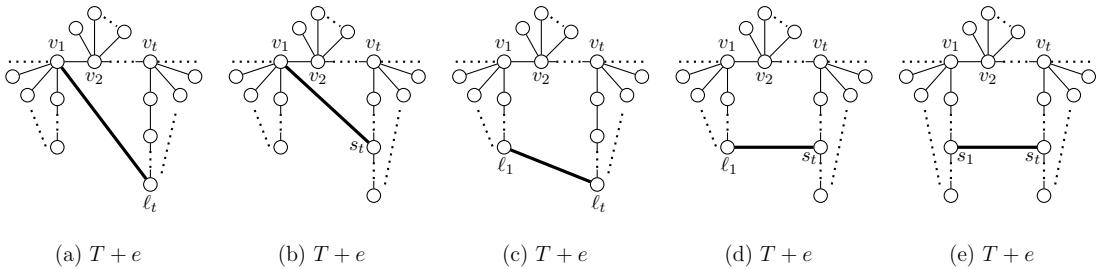

\emph{Subcase 4.2. $e=v_1\ell_t$ (see (a) of Figure~\ref{fig7}):} Take $e'=v_1v_2$. If $ter_T(v_t) \ge 3$, then $\ell_t$ becomes an interior degree two vertex in $T'$. If $ter_T(v_t)=2$, then $T'$ contains $v_t$ as an emv with $ter_{T'}(v_t)=1$.

\emph{Subcase 4.3. $e=v_1s_t$ (see (b) of Figure~\ref{fig7}):} If we take $e'=v_1v_2$, then $s_t$ becomes an emv with $ter_{T'}(s_t)=1$.

\emph{Subcase 4.4. $e=\ell_1\ell_t$ (see (c) of Figure~\ref{fig7}):} Let $e'=v_1v_2$. If $ter_T(v_1) \ge 3$ or $v_1$ is adjacent to at least two emvs in $T$, then $T'$ contains $\ell_1$ and $\ell_t$ as interior degree two vertices. If $ter_T(v_1)=2$ and $v_1$ is adjacent to exactly one emv $v_2$ in $T$, then $v_1$ becomes an exterior degree two vertex in $T'$ and no non-emv of $T$ becomes an emv in $T'$, and thus $ex(T')=k$.

\emph{Subcase 4.5. $e=\ell_1s_t$ (see (d) of Figure~\ref{fig7}):} If $ter_T(v_t) \ge 3$, take $e'=v_{t-1}v_t$; notice that $\ell_1$ is an interior degree two vertex and $s_t$ is an emv of terminal degree 1 in $T'$. (Similarly, if $ter_T(v_1) \ge 3$, take $e'=v_{1}v_2$; notice that $\ell_1$ is an interior degree two vertex and $s_t$ is an emv of terminal degree 1 in $T'$.) If $ter_T(v_1)=ter_T(v_t)=2$, then delete $e'=v_1x$, where $N(v_1) \cap V(C)=\{x, v_2\}$ (notice that $x$ could be $\ell_1$); then $T'$ has $v_t$ as an emv with $ter_{T'}(v_t)=1$.

\emph{Subcase 4.6. $e=s_1s_t$ (see (e) of Figure~\ref{fig7}):} If we take $e'=v_1v_2$, then $s_t$ becomes an emv with $ter_{T'}(s_t)=1$.\\

We have therefore proved that $\dim(T+e) \le Z(T+e)+1$ for any tree $T$ and for any $e \in E(\overline{T})$.~\hfill
\end{proof}

\begin{figure}[ht]
\centering
\begin{tikzpicture}[scale=.55, transform shape]
\node [draw, shape=circle] (1) at  (1,2) {};
\node [draw, shape=circle, fill=black] (2) at  (-0.3,1) {};
\node [draw, shape=circle] (3) at  (-0.3,0) {};
\node [draw, shape=circle, fill=black] (4) at  (-0.3,-1) {};
\node [draw, shape=circle, fill=black] (5) at  (1.5,1) {};
\node [draw, shape=circle] (6) at  (2.5,1) {};

\node [draw, shape=circle] (11) at  (7,2) {};
\node [draw, shape=circle, fill=black] (22) at  (5.7,1) {};
\node [draw, shape=circle] (33) at  (5.7,0) {};
\node [draw, shape=circle] (44) at  (5.7,-1) {};
\node [draw, shape=circle, fill=black] (55) at  (7.5,1) {};
\node [draw, shape=circle] (66) at  (8.5,1) {};

\node [scale=1.4] at (1,-1.5) {\large (a)};
\node [scale=1.4] at (7,-1.5) {\large (b)};

\draw(1)--(2)--(3)--(4)--(1);
\draw(5)--(1)--(6);
\draw(11)--(22)--(33)--(44)--(11);
\draw(55)--(11)--(66);

\end{tikzpicture}
\caption{A unicyclic graph $G$ with $\dim(G)=Z(G)+1$}\label{fig9}
\end{figure}
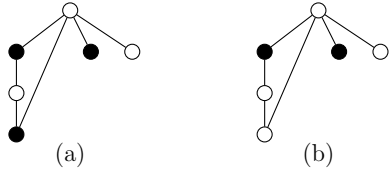

\begin{remark}
Figure \ref{fig9} shows a unicyclic graph $G$ with $\dim(G)=Z(G)+1$: The black vertices in (a) form a minimum resolving set of $G$, whereas the black vertices in (b) form a minimum zero forcing set of $G$.
\end{remark}

In analogy with the well-known notion of cycle rank, we define ``the even cycle rank of $G$", denoted by $r_e(G)$, to be the minimum number of edges to delete from $G$ so that the resulting graph $G'$ contains no even cycles. We originally conjectured that $\dim(G)\leq Z(G)+r_{e}(G)$, but after the submission of the first draft of this paper, we found the following counter-example.

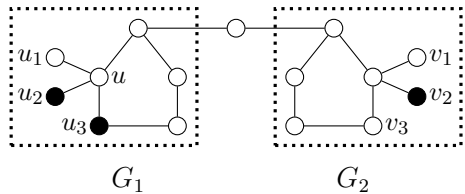
\begin{figure}[ht]
\centering
\begin{tikzpicture}[scale=.65, transform shape]
\node [draw, shape=circle] (1) at  (1,2) {};
\node [draw, shape=circle] (2) at  (0.2,1) {};
\node [draw, shape=circle, fill=black] (3) at  (0.2,0) {};
\node [draw, shape=circle] (4) at  (1.8,0) {};
\node [draw, shape=circle] (5) at  (1.8,1) {};
\node [draw, shape=circle] (a) at  (-0.7,1.4) {};
\node [draw, shape=circle, fill=black] (b) at  (-0.7,0.6) {};

\node [draw, shape=circle] (0) at  (3,2) {};

\node [draw, shape=circle] (11) at  (5,2) {};
\node [draw, shape=circle] (22) at  (4.2,1) {};
\node [draw, shape=circle] (33) at  (4.2,0) {};
\node [draw, shape=circle] (44) at  (5.8,0) {};
\node [draw, shape=circle] (55) at  (5.8,1) {};
\node [draw, shape=circle] (aa) at  (6.7,1.4) {};
\node [draw, shape=circle, fill=black] (bb) at  (6.7,0.6) {};

\node [scale=1.4] at (-1.2,1.4) {$u_1$};
\node [scale=1.4] at (-1.2,0.6) {$u_2$};
\node [scale=1.4] at (0.6,1) {$u$};
\node [scale=1.4] at (-0.3,0) {$u_3$};
\node [scale=1.4] at (7.2,1.4) {$v_1$};
\node [scale=1.4] at (7.2,0.6) {$v_2$};
\node [scale=1.4] at (6.3,0) {$v_3$};

\draw[very thick, style=dotted] (-1.6,-0.4) rectangle (2.2,2.4);
\draw[very thick, style=dotted] (3.8,-0.4) rectangle (7.6,2.4);
\node [scale=1.4] at (0.8,-1.1) {\large $G_1$};
\node [scale=1.4] at (5.4,-1.1) {\large $G_2$};

\draw(1)--(2)--(3)--(4)--(5)--(1);
\draw(1)--(0)--(11);
\draw(11)--(22)--(33)--(44)--(55)--(11);
\draw(a)--(2)--(b);
\draw(aa)--(55)--(bb);

\end{tikzpicture}
\caption{A graph $G$ with $r(G)=2$ and $r_e(G)=0$ such that $\dim(G)=4$ and $Z(G)=3$}\label{fig10}
\end{figure}

\begin{remark}\label{counter_ECRC}
There exists a graph $G$ containing no even cycles with $\dim(G)>Z(G)$ (see Figure \ref{fig10}); notice that $r(G)=2$ and $r_e(G)=0$. We show that $G$ in Figure \ref{fig10} satisfies $\dim(G)=4$ and $Z(G)=3$. Let $W$ and $S$ be the minimum resolving set and the minimum zero forcing set for $G$, respectively. First, we show that $\dim(G)=4$. Since $u_1$ and $u_2$ are twin vertices, $|W \cap \{u_1, u_2\}| \ge 1$, say $u_2 \in W$; similarly, we may assume that $v_2 \in W$. If $|W \cap V(G_1)|=1$, then $code_{W}(u_1)=code_W(u_3)$, and thus $|W \cap V(G_1)| \ge 2$; similarly, $|W \cap V(G_2)| \ge 2$. So, $|W| \ge 4$. Since $\{u_1,u_3,v_1,v_3\}$ forms a resolving set for $G$, we have $\dim(G)=4$. Next, we show that $Z(G)=3$. Notice that $S \cap \{u_1, u_2\} \neq \emptyset$; otherwise, the cut-vertex $u$ has two white neighbors $u_1$ and $u_2$, and thus $G$ fails to turn black after finitely many applications of the color-change rule. Similarly, $S \cap \{v_1, v_2\} \neq \emptyset$. WLOG, assume that $S_0=\{u_2, v_2\} \subseteq S$. Since $S_0$ fails to be a zero forcing set, $|S \setminus S_0| \ge 1$, and thus $Z(G) \ge 3$. Since $\{u_2,u_3,v_2\}$ forms a zero-forcing set, $Z(G) =3$.
\end{remark}

The forgoing example is particularly noteworthy because a further, more detailed analysis still (strongly) suggests that $\dim(G) \le Z(G)$ when the unique cycle of a unicyclic graph $G$ has odd length.

\begin{figure}[ht]
\centering
\begin{tikzpicture}[scale=.65, transform shape]
\node [draw, shape=circle] (1) at  (-0.5,2) {};
\node [draw, shape=circle] (2) at  (2,2) {};
\node [draw, shape=circle] (3) at  (4.5,2) {};
\node [draw, shape=circle] (4) at  (7,2) {};
\node [draw, shape=circle] (5) at  (9.5,2) {};
\node [draw, shape=circle] (6) at  (12,2) {};

\node [draw, shape=circle] (a1) at  (-0.5,1) {};
\node [draw, shape=circle] (a2) at  (-1,0.5) {};
\node [draw, shape=circle, fill=black] (a3) at  (-1,-0.2) {};
\node [draw, shape=circle] (a4) at  (0,-0.2) {};
\node [draw, shape=circle] (a5) at  (0,0.5) {};
\node [draw, shape=circle] (a6) at  (-1.6,0.8) {};
\node [draw, shape=circle, fill=black] (a7) at  (-1.6,0.2) {};

\node [draw, shape=circle] (b1) at  (2,1) {};
\node [draw, shape=circle] (b2) at  (1.5,0.5) {};
\node [draw, shape=circle] (b3) at  (1.5,-0.2) {};
\node [draw, shape=circle] (b4) at  (2.5,-0.2) {};
\node [draw, shape=circle] (b5) at  (2.5,0.5) {};
\node [draw, shape=circle] (b6) at  (0.9,0.8) {};
\node [draw, shape=circle, fill=black] (b7) at  (0.9,0.2) {};

\node [draw, shape=circle] (c1) at  (4.5,1) {};
\node [draw, shape=circle] (c2) at  (4,0.5) {};
\node [draw, shape=circle, fill=black] (c3) at  (4,-0.2) {};
\node [draw, shape=circle] (c4) at  (5,-0.2) {};
\node [draw, shape=circle] (c5) at  (5,0.5) {};
\node [draw, shape=circle] (c6) at  (3.4,0.8) {};
\node [draw, shape=circle, fill=black] (c7) at  (3.4,0.2) {};

\draw[style=solid] (4)--(7.4,2);
\draw[style=solid] (5)--(9,2);
\draw[style=dashed] (7.4,2)--(9,2);
\draw[style=dashed] (4)--(7,1);
\draw[style=dashed] (5)--(9.5,1);

\node [draw, shape=circle] (d1) at  (12,1) {};
\node [draw, shape=circle] (d2) at  (11.5,0.5) {};
\node [draw, shape=circle] (d3) at  (11.5,-0.2) {};
\node [draw, shape=circle] (d4) at  (12.5,-0.2) {};
\node [draw, shape=circle] (d5) at  (12.5,0.5) {};
\node [draw, shape=circle] (d6) at  (10.9,0.8) {};
\node [draw, shape=circle, fill=black] (d7) at  (10.9,0.2) {};

\draw[style=dashed] (-2,-0.5) rectangle (0.3,1.4);
\draw[style=dashed] (0.5,-0.5) rectangle (2.8,1.4);
\draw[style=dashed] (3,-0.5) rectangle (5.3,1.4);
\draw[style=dashed] (10.5,-0.5) rectangle (12.8,1.4);

\node [scale=1.4] at (-0.6,-0.9) {$G_1$};
\node [scale=1.4] at (1.9,-0.9) {$G_2$};
\node [scale=1.4] at (4.4,-0.9) {$G_3$};
\node [scale=1.4] at (11.9,-0.9) {$G_{2k}$};

\draw(1)--(2)--(3)--(4);
\draw(5)--(6);
\draw(1)--(a1)--(a2)--(a3)--(a4)--(a5)--(a1);
\draw(a6)--(a2)--(a7);
\draw(2)--(b1)--(b2)--(b3)--(b4)--(b5)--(b1);
\draw(b6)--(b2)--(b7);
\draw(3)--(c1)--(c2)--(c3)--(c4)--(c5)--(c1);
\draw(c6)--(c2)--(c7);
\draw(6)--(d1)--(d2)--(d3)--(d4)--(d5)--(d1);
\draw(d6)--(d2)--(d7);

\end{tikzpicture}
\caption{A graph $G$ with $r(G)=2k$ and $r_e(G)=0$ such that $\dim(G)=4k$ and $Z(G)=3k$}\label{fig11}
\end{figure}
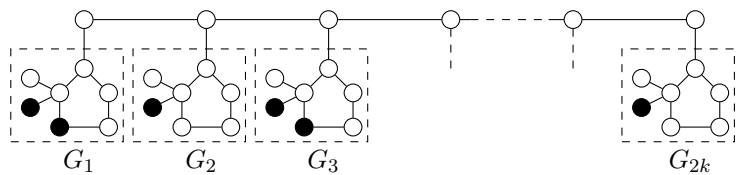

\begin{remark}
One can construct a graph $G$, containing no even cycles, such that $\dim(G)-Z(G)$ is arbitrary large (see Figure \ref{fig11}). Using the argument used in Remark \ref{counter_ECRC}, one can show that $\dim(G)=4k$ and $Z(G)=3k$ for the graph $G$ in Figure \ref{fig11}.
\end{remark}

Next, we recall a lemma which, together with results already shown, yields a more general inequality between metric dimension and zero forcing number.

\begin{lemma}[\cite{joc}]\label{joc_noeven}
If a connected graph $G$ has no even cycles and if $e$ is an edge of $G$, then $$\dim(G-e) \ge \dim(G)-1.$$
\end{lemma}

\begin{proposition}
Let $G$ be a connected graph containing no even cycles. Then $$\dim(G) \le Z(G)+2r(G).$$
\end{proposition}

\begin{proof}
Let $T$ be a spanning tree of $G$ obtained through the deletion of $r=r(G)$ edges of $G$. By Lemma~\ref{joc_noeven}, Theorem~\ref{DZtree}, and (b) of Theorem~\ref{Z+}, we have $\dim(G)\le \dim(T)+r \le Z(T)+r \le Z(G)+2r$.~\hfill
\end{proof}

Though without a proof, we are inclined to think that a stronger inequality exists between $\dim(G)$ and $Z(G)$ in full generality.

\begin{conjecture}\label{ECRConjecture}
(Cycle Rank Conjecture): $\dim(G)\leq Z(G)+r(G)$.
\end{conjecture}

\begin{figure}[ht]
\centering
\begin{tikzpicture}[scale=.65, transform shape]
\node [draw, shape=circle, fill=black] (1) at  (-1.8,0) {};
\node [draw, shape=circle] (2) at  (0,0) {};
\node [draw, shape=circle] (3) at  (1.8,0) {};
\node [draw, shape=circle, fill=black] (a1) at  (-1.1,-0.5) {};
\node [draw, shape=circle] (a2) at  (-1.4,-1.2) {};
\node [draw, shape=circle] (a3) at  (-0.5,-1.1) {};
\node [draw, shape=circle] (b1) at  (1.1,-0.5) {};
\node [draw, shape=circle] (b2) at  (1.4,-1.2) {};
\node [draw, shape=circle, fill=black] (b3) at  (0.5,-1.1) {};
\node [draw, shape=circle, fill=black] (c1) at  (1.1,0.5) {};
\node [draw, shape=circle] (c2) at  (1.4,1.2) {};
\node [draw, shape=circle] (c3) at  (0.5,1.1) {};
\draw[very thick, style=dotted] (-0.9, -0.1).. controls (-0.5,0.6).. (0.2, 0.8);

\node [draw, shape=circle, fill=black] (11) at  (5.2,0) {};
\node [draw, shape=circle] (22) at  (7,0) {};
\node [draw, shape=circle] (33) at  (8.8,0) {};
\node [draw, shape=circle, fill=black] (a11) at  (5.9,-0.5) {};
\node [draw, shape=circle] (a22) at  (5.6,-1.2) {};
\node [draw, shape=circle, fill=black] (a33) at  (6.5,-1.1) {};
\node [draw, shape=circle, fill=black] (b11) at  (8.1,-0.5) {};
\node [draw, shape=circle] (b22) at  (8.4,-1.2) {};
\node [draw, shape=circle, fill=black] (b33) at  (7.5,-1.1) {};
\node [draw, shape=circle, fill=black] (c11) at  (8.1,0.5) {};
\node [draw, shape=circle] (c22) at  (8.4,1.2) {};
\node [draw, shape=circle, fill=black] (c33) at  (7.5,1.1) {};
\draw[very thick, style=dotted] (6.1, -0.1).. controls (6.5,0.6).. (7.2, 0.8);

\node [scale=1.4] at (0,-2.2) {\large $Z(G)=k+1$};
\node [scale=1.4] at (7,-2.2) {\large $\dim(G)=2k+1$};

\draw(1)--(2)--(3);
\draw(2)--(a1)--(a2)--(a3)--(2);
\draw(2)--(b1)--(b2)--(b3)--(2);
\draw(2)--(c1)--(c2)--(c3)--(2);
\draw(11)--(22)--(33);
\draw(22)--(a11)--(a22)--(a33)--(22);
\draw(22)--(b11)--(b22)--(b33)--(22);
\draw(22)--(c11)--(c22)--(c33)--(22);

\end{tikzpicture}
\caption{A graph $G$ satisfying $\dim(G)=Z(G)+r(G)$}\label{fig12}
\end{figure}
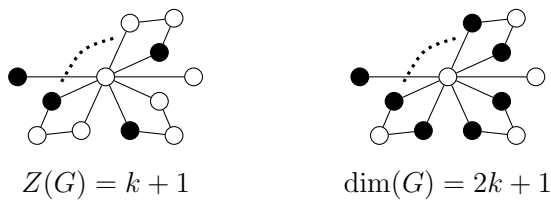

\begin{remark}
There exists a graph $G$ satisfying $\dim(G)=Z(G)+r(G)$ (see Figure~\ref{fig12}). Suppose that $G$ is a graph obtained by identifying $k$ copies of $C_4$ to the central vertex of $P_3$. Then one can readily verify that $Z(G)=k+1$ and $\dim(G)=2k+1$; also notice that $r(G)=k$. Thus, we have $\dim(G)=Z(G)+r(G)$.
\end{remark}

\section{A Proof of Theorem~\ref{unicyclic}}

In this final section, for reasons explained in the second paragraph of Section 3, we include a proof of Theorem~\ref{unicyclic}. We begin with counter-examples to two important assertions made in the outline of proof in~\cite{CEJO} to Theorem~\ref{unicyclic}.

\begin{figure}[ht]
\centering
\begin{tikzpicture}[scale=.65, transform shape]
\node [draw, shape=circle] (1) at  (0,0.6) {};
\node [draw, shape=circle] (2) at  (1,0) {};
\node [draw, shape=circle] (3) at  (2,0) {};
\node [draw, shape=circle] (4) at  (3,0) {};
\node [draw, shape=circle] (5) at  (0,-0.6) {};
\node [draw, shape=circle] (6) at  (-1,-1.2) {};

\node [draw, shape=circle] (11) at  (6,0.6) {};
\node [draw, shape=circle] (22) at  (7,0) {};
\node [draw, shape=circle] (33) at  (8,0) {};
\node [draw, shape=circle] (44) at  (9,0) {};
\node [draw, shape=circle] (55) at  (6,-0.6) {};
\node [draw, shape=circle] (66) at  (5,-1.2) {};

\node [scale=1.4] at (7.1,-.9) {$e$};
\node [scale=1.4] at (1,-2.2) {\large $T$};
\node [scale=1.4] at (7,-2.2) {\large $T+e$};

\draw(1)--(2)--(3)--(4);
\draw(6)--(5)--(2);
\draw(11)--(22)--(33)--(44);
\draw(66)--(55)--(22);
\draw[very thick, black] (55).. controls (7.2,-0.55).. (33);

\end{tikzpicture}
\caption{Unicyclic graph $T+e$ such that $ex(T+e)>ex(T)$}\label{fig13}
\end{figure}
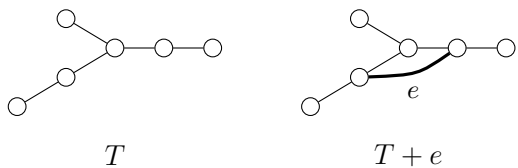

\begin{remark} \label{ce1CEJO} Figure~\ref{fig13} shows a counter-example to the assertion ``$ex(T+e)\leq ex(T)$" in the first sentence of the outline of proof: note that $ex(T)=1$ and $ex(T+e)=3$.
\end{remark}

\begin{figure}[ht]
\centering
\begin{tikzpicture}[scale=.65, transform shape]
\node [draw, shape=circle] (1) at  (-1.5,1) {};
\node [draw, shape=circle] (2) at  (0,1) {};
\node [draw, shape=circle] (3) at  (1.5,1) {};
\node [draw, shape=circle] (4) at  (3,1) {};
\node [draw, shape=circle] (5) at  (4.5,1) {};
\node [draw, shape=circle] (6) at  (6,1) {};
\node [draw, shape=circle] (7) at  (7.5,1) {};
\node [draw, shape=circle] (8) at  (9,1) {};

\node [draw, shape=circle] (a1) at  (0,-0.2) {};
\node [draw, shape=circle] (a2) at  (0,-1.4) {};
\node [draw, shape=circle] (b1) at  (1.5,-0.2) {};
\node [draw, shape=circle] (b2) at  (1.5,-1.4) {};
\node [draw, shape=circle] (e1) at  (6,-0.2) {};
\node [draw, shape=circle] (e2) at  (6,-1.4) {};
\node [draw, shape=circle] (f1) at  (7.5,-0.2) {};
\node [draw, shape=circle] (f2) at  (7.5,-1.4) {};
\node [draw, shape=circle] (c1) at  (3,2.2) {};
\node [draw, shape=circle] (c2) at  (3,3.4) {};
\node [draw, shape=circle] (d1) at  (4.5,2.2) {};
\node [draw, shape=circle] (d2) at  (4.5,3.4) {};

\node [scale=1.4] at (3.7,-0.5) {\large $e$};

\draw(1)--(2)--(3)--(4)--(5)--(6)--(7)--(8);
\draw(2)--(a1)--(a2);
\draw(3)--(b1)--(b2);
\draw(6)--(e1)--(e2);
\draw(7)--(f1)--(f2);
\draw(4)--(c1)--(c2);
\draw(5)--(d1)--(d2);
\draw[very thick](b1)--(e1);

\end{tikzpicture}
\caption{Unicyclic graph $T+e$ satisfying the assumption of Case 4 in the outline of proof such that $W \neq \emptyset$ (see p.109 of~\cite{CEJO} for the definition of $W$) and $T$ is not a caterpillar}\label{fig14}
\end{figure}
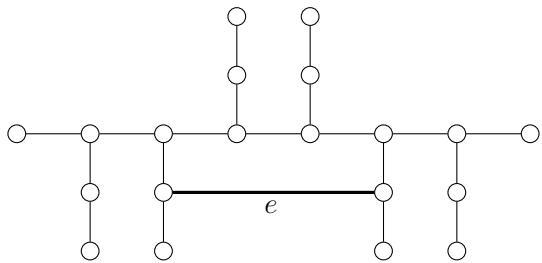

\begin{remark} \label{ce2CEJO} Figure~\ref{fig14} shows a counter-example to the argument given in the outline of proof (Case 4).
\end{remark}

\textit{Proof of Theorem~\ref{unicyclic}.}\\

We first prove the lower bound; namely, $\dim(T)-2\leq \dim(T+e)$ where $e\in E(\overline{T})$. Since the inequality obviously holds for a path, let $T$ be a tree which is not a path. Recall that $\dim(T+e)\geq \sigma(T+e)-ex(T+e)$ by Lemma~\ref{dim geq sigma-ex} and $\dim(T)=\sigma(T)-ex(T)$ by Theorem~\ref{tree}. Let $e=uv$; we consider three cases pending the nature of vertices $u$ and $v$. Case 1: $u$ and $v$ are both end-vertices. In this case, $\sigma(T+e)=\sigma(T)-2$ (sum of terminal degrees reduces by two) and $ex(T+e)\leq ex(T)$ (no new exterior major vertices are created). Thus, $\dim(T+e)\geq \sigma(T+e)-ex(T+e) \geq \sigma(T)-2-ex(T)=\dim(T)-2$. Case 2: exactly one of $u$ and $v$ is an end-vertex. In this case, $\sigma(T+e)=\sigma(T)-1$ and $ex(T+e)\leq ex(T)+1$. Thus, $\dim(T+e)\geq \sigma(T+e)-ex(T+e) \geq \sigma(T)-1-(ex(T)+1)=\dim(T)-2$. Case 3: neither $u$ nor $v$ is an end-vertex. In this case, $\sigma(T+e)=\sigma(T)$ and $ex(T+e)\leq ex(T)+2$. Thus, $\dim(T+e)\geq \sigma(T+e)-ex(T+e) \geq \sigma(T)-(ex(T)+2)=\dim(T)-2$. \\

Now, we prove the upper bound; namely, $\dim(T+e)\leq \dim(T)+1$ where $e\in E(\overline{T})$. The claim holds when $T$ is a path $P_n$, as the two end-vertices of $P_n$ form a basis (minimum resolving set) for $P_n+e$: If $e=v_iv_j$ where $i<j$, then $v_i$ and $v_j$, being adjacent vertices, resolve vertices on the unique cycle $C$ of $P_n+e$ among themselves (whence we say ``$v_i$ and $v_j$ resolve $C$"). But then $W=\{v_1,v_n\}$ resolves $C$ since for any $v \in V(C)$, $code_{W'}(v)=code_W(v)+(a_1, a_2)$, where $W'=\{v_i, v_j\}$ and $(a_1, a_2)$ is a fixed vector. Further, $v_1$ and $v_n$ obviously resolve vertices in $V(P_n+e)\setminus V(C)$ among themselves and from $V(C)$.\\

So, let $T$ be a tree which is not a path, and thus $\dim(T)\geq 2$. Cyclically label the vertices lying on the unique cycle $C$ of $T+e$ ($e\in E(\overline {T})$) by $u_1, \ldots, u_k$ ($k\geq 3$). Denote by $T_i$ the subtree rooted at $u_i$ (in other words, the component of $(T+e)\setminus E(C)$ which contains $u_i$). Given any basis $B$ of $T$, partition $B$ into the disjoint union of sub-bases $B_i$, where $B_i\subseteq V(T_i)$, $1\leq i\leq k$; assume, without loss of generality, that $B_1\neq \emptyset$. If $B_i=\emptyset$ for each $i\neq 1$, then $T- T_1$ must be a path (for $B$ to be a basis of $T$); in this case, either $B\cup\{u_2\}$ or $B\cup\{u_k\}$ is a resolving set for $T+e$.\\

So, assume there exists $1< i \leq k$ such that $B_i\neq \emptyset$. If there exist two non-empty sub-bases $B_i$ and $B_j$ such that $d_{T+e}(u_i,u_j)=m=\lfloor\frac{k}{2}\rfloor$, then let $b_0 \in V(C)\setminus \{u_i,u_j\}$ and put $B_0=\{b_i,b_j,b_0\}$ (also put $B'_0=\{u_i,u_j,b_0\}$) where $b_i\in B_i$ and $b_j\in B_j$; otherwise, let $b_0=u_{m+1}$ and put $B_0=\{b_1, b_0, b_s\}$ (also put $B'_0=\{u_1,b_0,u_s\}$), where $b_1\in B_1$ and $b_s\in B_s\neq \emptyset$ for some $s\neq 1,m+1$. (The point here is to arrange a resolving set for $T+e$ that contains elements in three subtrees (the $T_i$'s), two of which having roots (the $u_i$'s) attaining the diameter of the cycle $C$.) We will show that the set $\widetilde{B}=B\cup\{b_0\}$ is a resolving set for $T+e$. Notice that $B_0\subseteq \widetilde{B}$.\\

By Lemma~\ref{resolving subtrees}, we have $code_{B_0}(x_i)\neq code_{B_0}(x_j)$ and, a fortiori, $code_{\widetilde{B}}(x_i)\neq code_{\widetilde{B}}(x_j)$ for $x_i\in V(T_i)$ and $x_j\in V(T_j)$, when $i\neq j$. It thus suffices to show that $\forall x,y \in V(T_i)$ where $1\leq i \leq k$, $code_{\widetilde{B}}(x)\neq code_{\widetilde{B}}(y)$. Accordingly, let $x, y\in V(T_i)$ be given for a fixed $i$. It's clear that if $d_{T}(x,b)\neq d_{T}(y,b)$ for some $b\in B_i$, then $d_{T+e}(x,b)\neq d_{T+e}(y,b)$; so, let $b\in B_j$ for some $j\neq i$. Notice that there exists a fixed $a\in \mathbb{N}$ such that $\forall x\in V(T_i)$, $d_{T+e}(x,b)=d_{T}(x,b)-a$. Thus, $d_{T}(x,b)\neq d_{T}(y,b)$ implies $d_{T+e}(x,b)\neq d_{T+e}(y,b)$ for $b\notin B_i$ as well.\\

We have thus proved the theorem.~\hfill$\Box$\\

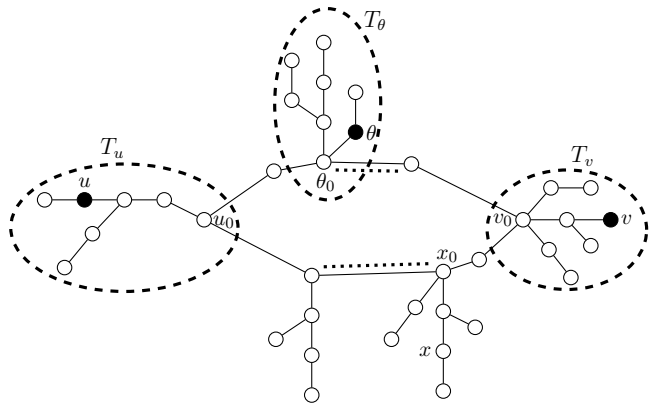
\begin{figure}[ht]
\centering
\begin{tikzpicture}[scale=.53, transform shape]

\draw[very thick, style=dashed](-0.9,2.9) ellipse (1.33 and 2.4);
\draw[very thick, style=dashed](-6,-0.2) ellipse (2.85 and 1.6);
\draw[very thick, style=dashed](5.1,-0.25) ellipse (2 and 1.5);

\node [draw, shape=circle] (a0) at  (-1,1.45) {};
\node [draw, shape=circle] (1) at  (-2.25,1.22) {};
\node [draw, shape=circle] (b0) at  (-4,0) {};
\node [draw, shape=circle] (c0) at  (-1.3,-1.4) {};
\node [draw, shape=circle] (d0) at  (2,-1.3) {};
\node [draw, shape=circle] (2) at  (2.9,-1) {};
\node [draw, shape=circle] (e0) at  (4,0) {};
\node [draw, shape=circle] (3) at  (1.2,1.4) {};

\draw[very thick, style=dotted](-1, -1.18)--(1.7, -1.1); 
\draw[very thick, style=dotted](-0.7, 1.25)--(0.9, 1.22); 

\node [draw, shape=circle] (a1) at  (-1,2.45) {};
\node [draw, shape=circle] (a2) at  (-1,3.45) {};
\node [draw, shape=circle] (a3) at  (-1,4.45) {};
\node [draw, shape=circle] (a11) at  (-1.8,3) {};
\node [draw, shape=circle] (a22) at  (-1.8,4) {};
\node [draw, shape=circle, fill=black] (a111) at  (-0.2,2.2) {};
\node [draw, shape=circle] (a222) at  (-0.2,3.2) {};

\node [draw, shape=circle] (b1) at  (-5,0.5) {};
\node [draw, shape=circle] (b2) at  (-6, 0.5) {};
\node [draw, shape=circle, fill=black] (b3) at  (-7,0.5) {};
\node [draw, shape=circle] (b4) at  (-8,0.5) {};
\node [draw, shape=circle] (b11) at  (-6.8,-0.35) {};
\node [draw, shape=circle] (b22) at  (-7.5,-1.2) {};

\node [draw, shape=circle] (c1) at  (-1.3,-2.4) {};
\node [draw, shape=circle] (c2) at  (-1.3, -3.4) {};
\node [draw, shape=circle] (c3) at  (-1.3,-4.4) {};
\node [draw, shape=circle] (c11) at  (-2.2,-3) {};

\node [draw, shape=circle] (d1) at  (2,-2.3) {};
\node [draw, shape=circle] (d2) at  (2,-3.3) {};
\node [draw, shape=circle] (d3) at  (2,-4.3) {};
\node [draw, shape=circle] (d11) at  (1.3,-2.2) {};
\node [draw, shape=circle] (d22) at  (0.7,-3) {};
\node [draw, shape=circle] (d111) at  (2.8,-2.7) {};

\node [draw, shape=circle] (e1) at  (5.1,0) {};
\node [draw, shape=circle, fill=black] (e2) at  (6.2,0) {};
\node [draw, shape=circle] (e11) at  (4.65,-0.75) {};
\node [draw, shape=circle] (e22) at  (5.2,-1.45) {};
\node [draw, shape=circle] (e111) at  (5.7,-0.65) {};
\node [draw, shape=circle] (e1111) at  (4.7,0.8) {};
\node [draw, shape=circle] (e2222) at  (5.7,0.8) {};

\node [scale=1.4] at (-0.95,0.95) {$\theta_0$};
\node [scale=1.4] at (0.2,2.2) {$\theta$};
\node [scale=1.4] at (-3.5,-0.05) {$u_0$};
\node [scale=1.4] at (-7,0.95) {$u$};
\node [scale=1.4] at (2.1,-0.9) {$x_0$};
\node [scale=1.4] at (1.55,-3.3) {$x$};
\node [scale=1.4] at (3.5,0) {$v_0$};
\node [scale=1.4] at (6.6,0) {$v$};

\node [scale=1.4] at (0.3,5) {\large $T_{\theta}$};
\node [scale=1.4] at (-6.3,1.8) {\large $T_u$};
\node [scale=1.4] at (5.5,1.6) {\large $T_v$};

\draw(a3)--(a2)--(a1)--(a0)--(a111)--(a222);
\draw(a1)--(a11)--(a22);
\draw(b0)--(b1)--(b2)--(b3)--(b4);
\draw(b2)--(b11)--(b22);
\draw(c0)--(c1)--(c2)--(c3);
\draw(c1)--(c11);
\draw(d22)--(d11)--(d0)--(d1)--(d2)--(d2)--(d3);
\draw(d1)--(d111);
\draw(e22)--(e11)--(e0)--(e1)--(e2);
\draw(e1)--(e111);
\draw(e0)--(e1111)--(e2222);

\draw(a0)--(1)--(b0)--(c0)--(d0)--(2)--(e0)--(3)--(a0);

\end{tikzpicture}
\caption{The set $\{u,v,\theta\}$ resolves the subtrees $T_i$'s from each other}\label{fig15}
\end{figure}

The following lemma shows that subtrees are distinguished by the $B_0$ chosen above; see Figure~\ref{fig15} for an illustration of the situation under consideration.

\begin{lemma}\label{resolving subtrees}
Let $B_0$ and $B'_0$ be chosen as in the Proof of Theorem~\ref{unicyclic}; explicitly, let $B_0=\{u,v,\theta\}$ and $B'_0=\{u_0,v_0,\theta_0\}\subseteq V(C)$, where $d(u_0,v_0)=diam(C)$ and $u$ ($v, \theta$, respectively) is a vertex on the subtree rooted at $u_0$ ($v_0, \theta_0$, respectively). Then, we have $code_{B_0}(x)\neq code_{B_0}(y)$ for vertices $x$ and $y$ belonging to distinct subtrees rooted at vertices of the unique cycle $C$ of $T+e$.
\end{lemma}

\begin{proof}
Observe that $B'_0$ \emph{strongly resolves} the unique cycle $C$ of $T+e$ (see the first paragraph of Section 3 for the definition of ``strongly resolves"), because no vertex of $C$ can have shorter distance, \emph{by the same value}, to all vertices of $B_0'$ than another vertex of $C$. Thus, $B_0$ \emph{strongly resolves} $C$, because there exists a fixed vector $(a_1,a_2,a_3)$ such that $\forall x\in V(C)$, $code_{B_0}(x)=code_{B_0'}(x)+ (a_1,a_2,a_3)$. If $x\in V(T_i)$ where $V(T_i)\cap B_0=\emptyset$, then $[x]_{B_0}=[x_0]_{B_0}$, where $x_0$ is the root of $T_i$: this is because any path from $x$ \emph{of such a subtree $T_i$} to a vertex in $B_0$ must go through $x_0$. Thus $[x]_{B_0}\neq [y]_{B_0}$ and, a fortiori, $code_{B_0}(x)\neq code_{B_0}(y)$ for $x$ and $y$ belonging to distinct subtrees which have empty intersection with $B_0$. If $B_0=B'_0$, then the same reasoning applies to the subtrees containing elements of $B_0$. Otherwise, if suffices to check $code_{B_0}(x)\neq code_{B_0}(y)$ (1) for $x\in V(T_i)$ and $y\in V(T_u)$, (2) for $x\in V(T_i)$ and $y\in V(T_\theta)$, (3) for $x\in V(T_u)$ and $y\in V(T_v)$, and (4) for $x\in V(T_u)$ and $y\in V(T_\theta)$; here $T_u, T_v$, $T_\theta$, and $T_i$ are the subtrees containing $u, v$, $\theta$, and none of $B_0$, respectively. Since the same argument works for all four inequalities, we will only explicitly verify (1). \\

Suppose, for the sake of contradiction, $code_{B_0}(y)=code_{B_0}(x)$; i.e., $(d(y,u),d(y,v),d(y,\theta))=(d(x,u),d(x,v),d(x,\theta))$ for vertices $y\in V(T_u)$ and $x\in V(T_i)$. Equating the first two coordinates and expanding, we get $d(y,u)=d(x,x_0)+d(x_0,u_0)+d(u_0,u)$ and $d(y,u_0)+d(u_0,v_0)+d(v_0,v)=d(x,x_0)+d(x_0,v_0)+d(v_0,v)$, where $x_0$ is the root of the subtree containing $x$. Subtracting the two equations and rearranging terms, we get $d(y,u)=d(y,u_0)+d(x_0,u_0)+d(u_0,u)+d(u_0,v_0)-d(x_0,v_0)$. Now, since $d(u_0,v_0)=diam(C)$, we have $d(u_0,v_0)-d(x_0,v_0)=d(u_0,x_0)$. And we have $d(y,u)=d(y,u_0)+d(u_0,u)+2d(u_0,x_0)$. Since $x\in V(T_i)$ and $T_i\neq T_u$, $d(u_0,x_0) > 0$, and we have $d(y,u) > d(y,u_0)+d(u_0,u)$, violating the triangle inequality which $d(\cdot,\cdot)$ must satisfy as a metric.~\hfill
\end{proof}

\begin{remark}
Notice that Lemma~\ref{resolving subtrees} still holds if each ``subtree $T_i$ rooted at $u_i$" is replaced by ``subgraph $G_i$ rooted at $u_i$" with $G_i$ and $G_j$ disjoint for $i \neq j$.
\end{remark}


\noindent\textbf{Acknowledgements} The authors wish to thank the anonymous referees for their comments and suggestions.

\end{document}